\documentclass[10pt]{article}
\usepackage{fullpage,epsfig,graphics,amsbsy,amssymb,cancel}
\usepackage{psfrag,hyperref}
\usepackage{graphicx}
\usepackage{tikz}
\usepackage{amsfonts}
\usepackage{amssymb,amsmath,amscd}

\newcommand{\BS}{\boldsymbol}
\newcommand{\BB}{\mathbb}

\newcommand{\bea}{\begin{eqnarray}}
\newcommand{\eea}{\end{eqnarray}}
\newcommand{\beq}{\begin{eqnarray}}
\newcommand{\eeq}{\end{eqnarray}}
\newcommand{\nn}{\nonumber}

\def\L{{\cal L}}
\def\M{{\cal M}}
\def\N{{\cal N}}
\def\R{{\cal R}}

\def\E{{\cal E}}

\def\So{{\cal X}}
\def\Ta{{\cal M}}
\def\Pa{{\cal L}}
\def\Mod{\textsf{M}}
\def\Nod{\textsf{N}}
\def\A{\textsf{A}}

\newcommand{\smalint}{{\Large\textrm{$\smallint$}}}
\newcommand{\sbullet}{\textrm{\tiny{\textbullet}}}

\newcommand{\FR}{\mathfrak}
\def\ga{\alpha}
\def\gb{\beta}
\def\gc{\gamma}
\def\gd{\delta}
\def\ge{\epsilon}
\newcommand{\Hom}{\textrm{Hom}}
\newcommand{\rep}{{\cal R}\textrm{ep}^{\infty}}
\newcommand{\Set}{\textsf{Set}}
\newcommand{\HomM}{\widehat{\textrm{Hom}}}

\newcommand{\cdga}{\textsf{cdga}}
\newcommand{\dga}{\textsf{dga}}
\newcommand{\gam}{\textsf{cga}}
\newcommand{\dgm}{\textsf{dgM}}

\newtheorem{Thm}{Theorem}[section]
\newtheorem{Prop}[Thm]{Proposition}
\newtheorem{Lem}[Thm]{Lemma}
\newtheorem{Cor}[Thm]{Corollary}

\newtheorem{Exa}[Thm]{Example}

\newtheorem{lem}[Thm]{Lemma}
\newtheorem{cor}[Thm]{Corollary}

\newtheorem{rmk}[Thm]{Remark}
\newtheorem{defi}[Thm]{Definition}
\usetikzlibrary{arrows,chains,matrix,positioning,scopes}
\makeatletter
\tikzset{join/.code=\tikzset{after node path={%
\ifx\tikzchainprevious\pgfutil@empty\else(\tikzchainprevious)%
edge[every join]#1(\tikzchaincurrent)\fi}}}
\makeatother

\tikzset{>=stealth',every on chain/.append style={join},
         every join/.style={->}}

\begin{document}

\thispagestyle{empty}
\begin{flushright} \small
UUITP-22/11\\
NSF-KITP-11-170\\
 \end{flushright}
\smallskip
\begin{center} \Large
{\bf Wilson Lines from Representations of NQ-Manifolds}
  \\[12mm] \normalsize
{\bf Francesco Bonechi$^a$}, {\bf Jian Qiu$^a$ and Maxim Zabzine$^{b,c}$} \\[8mm]
 {\small\it
 ${}^a$I.N.F.N. and Dipartimento di Fisica\\
     Via G. Sansone 1, 50019 Sesto Fiorentino - Firenze, Italy\\
   \vspace{.3cm}
     ${}^b$Department of Physics and Astronomy,
     Uppsala university,\\
     Box 516,
     SE-751\;20 Uppsala,
     Sweden\\
      \vspace{.3cm}
     ${}^c$Kavli Institute for Theoretical Physics, University of California, \\
Santa Barbara, CA 93106 USA\\
 }
\end{center}
\vspace{7mm}

\begin{abstract}
\noindent An $NQ$-manifold is a non-negatively graded supermanifold with a degree 1 homological vector field. The focus of this paper is to define the Wilson loops/lines in the context of $NQ$-manifolds and to study their properties.
The Wilson loops/lines, which give the holonomy or parallel transport, are familiar objects in usual differential geometry, we analyze the subtleties in
the generalization to the $NQ$-setting and we also sketch some possible applications of our construction.
\end{abstract}
\tableofcontents

\section{Introduction}\label{Intro}
An $NQ$-manifold $(\Ta,Q)$ is a non-negatively graded supermanifold $\Ta$ with a degree 1 homological vector field $Q$.
A representation of a $NQ$ manifold $\cal M$ is a graded vector bundle $\E$ over ${\cal M}$ endowed with a lift of the $Q$-structure linear in the
fibre coordinates
(a planar $Q$ connection in \cite{vaintrob}). As an example, consider ${\cal M}=L[1]$ to be a Lie algebroid, and $\E=T^*[2]L[1]$ endowed with
the cotangent lift of $Q$; this construction gives a canonical formulation of the
adjoint representation
which is described as a representation up to homotopy in \cite{abad-2009}.
The major goal of this paper is to lay down the notion of Wilson loops and lines for representations of $NQ$-manifolds.

We always resort to the intuition given by the standard differential geometry description of flat connections, where the Wilson loops and
Wilson lines measure holonomy and parallel transport. On the other hand, when the transition functions of $\E$ have components in degree
higher than zero, issues such as the trivialization (in)dependence\footnote{{\color{black}By trivialization independence, it is meant that the Wilson-line transforms
covariantly under a change of the trivialization, so that one may glue pieces of Wilson lines together in a consistent manner.}}
and homotopy (in)variance become very delicate.


A motivation for this paper is to clarify certain issues left over from the work done by two of the authors in ref.\cite{2010arXiv1006.1240Q}.
We first constructed Wilson-loops for an $NQ$-manifold and showed that by evaluating the expectation value of the Wilson-loop,
one obtain new weight systems for knots embedded in 3-manifolds. But the issue of trivialization dependence was not discussed
thoroughly and we feel the need to write a separate paper now and examine these issues with painstaking care as the concept of
Wilson-lines/loops for graded geometry is rather new and unfamiliar.

Furthermore, our construction may also have some bearing on the integration problem.
The main idea about the integration of $NQ$ manifolds goes back to ref.\cite{Severa01}.
A presheaf is defined as a contravariant functor from a category ${\cal C}$ to the category of sets $\Set$.
For an $NQ$-manifold ${\cal M}$ we can consider the following presheaf (which is the adjoint of the differentiation functor 1-jet of \cite{Severa06})
on the category of smooth manifolds Man:
\bea \smalint \Ta:~\textrm{Man}\to\Set,~~~~\smalint{\Ta}(\So)=\Hom_{NQ}(T[1]\So, \Ta),\label{significant}\eea
where $\So$ is a smooth manifold and $\Hom_{NQ}(T[1]\So, \Ta)$ is  the set of $NQ$-morphisms $T[1]\So\to\Ta$; this set is called the \emph{Maurer-Cartan set}.
Picking the standard simplices $\Delta^n$, one may describe $\int \Ta$ as a simplicial manifold
\bea (\smalint \Ta)_n:~~[n]\to \Hom_{NQ}(T[1]\Delta^n,\Ta),\label{def_MC_simp}\eea
$(\int\Ta)_{\sbullet}$ can be thought of as the singular chains of $\int \Ta$. In the case ${\cal M}=L[1]$, $(\int L[1])_{\sbullet}$ is called the
$\infty$-groupoid of $L$ and denoted as $\Pi_\infty(L)$; also recall that $\pi_1(L)=\Pi_1(L)/({\rm homotopy})$ is the $1$-groupoid of $L$ \cite{Severa01}.
In ref.\cite{AbadSchaetz} the integration of representations up to homotopy of $L$ into representations of $\Pi_\infty(L)$ is studied;
their paper and the present one have quite an overlap and yet are from decidedly different perspectives.
When $\Ta$ comes from an $L_{\infty}$-algebra,
various properties of $(\int\Ta)_{\sbullet}$ have been studied by Hinich \cite{1996alg.geom.6010H}, Getzler \cite{Getzler} and
Henriques \cite{Henriques}, in particular the Maurer-Cartan set is a Kan complex, which one recalls is the prominent trait of
the singular complex of a topological space.

In this paper, we choose an alternative approach to the simplicial set route which we find more natural in the $NQ$-manifold setting.
We consider in fact the presheaf on $NQ$
\bea \HomM_{NQ}(T[1][0,1],\Ta): NQ\rightarrow \Set,~~~ {\cal L}\rightarrow \Hom_{NQ}(T[1][0,1]\times {\cal L},\Ta).\label{OURpresheaf}\eea
In short, we look at $NQ$-maps from $T[1][0,1]\times \Pa$ to $\M$, namely super-paths 'parameterized' by a test manifold $\Pa$.
We can think that, by specializing the test manifolds ${\cal L}$ to the $n$-cubes, this presheaf defines the singular cubical complex of $\Ta$.
It is shown in prop.\ref{representability} that the presheaf is representable by an $NQ$-manifold $(\bar {\cal S},\delta_B)$.
Moreover, we subscribe to the the approach of Fock and Rosly \cite{FockRosly} along with Andersen et
al \cite{Andersen} who described the moduli space of flat connections over Riemann surfaces with punctures by the algebraic relations between
the Wilson loops. Indeed, they showed that for certain gauge groups,
the Wilson loops with possible double points (called chord diagrams) form a complete basis
of the algebraic functions on the moduli space of flat connections.

So we look at Wilson loops as functions on $\bar{\cal S}_{\textrm{O}}$ (which is similar to $\bar{\cal S}$ above but for closed loops) and at Wilson lines as parallel transport
 of $\E$ depending on $\bar{\cal S}$
(see Section \ref{sec_WL} for definitions).
Our construction of the Wilson-line is largely based on Chen's iterated integral used widely in the literature, yet it differs from the others' approach
quite significantly in the choice of the presheaf Eq.\ref{OURpresheaf}. Our main result is that Wilson loops are invariant under the change of trivialization
so that they define degree zero functions on $\bar{\cal S}_{\textrm{O}}$; moreover they are homotopy invariant.
In contrast, Wilson lines are not covariant under the change of trivialization; yet can be made so if we restrict the test manifold $\Pa$ in Eq.\ref{OURpresheaf} above to have zero $Q$-structure.
The restriction says that the Wilson-line only gives a presheaf on a subcategory of $NQ$. In many aspects, this restriction discards vast amount of information of the space of
super-curves, whilst some other less restrictive options do exist but will be on the case by case basis and will not be given detailed discussion here.

The paper is organized as follows:
we first review the notion of the representation of an $NQ$-manifold given by Vaintrob in sec.\ref{RoQS}, from these representations,
we can form the parallel transport operator using iterated integral. But to handle the graded objects properly, we spend some pages to review
the functor of points view of the graded manifolds, whose embodiment is the superfields that are widely used in physics computations.
The discussion in sec.\ref{BCaII} about the bar resolution is a bit of a digression, but has a bearing upon the problem: are the Wilson loops a
sufficient basis for the functions of the generalized connections. In the quite weighty sec.\ref{sec_WL}, we investigate the problem of trivialization
independence and homotopy invariance and point out some potential problems regarding these issues in the literature. In sec.\ref{sec_ingetnamn},
we sketch how one may apply the Wilson-loops in the integration problem, and possible interpretations of the Wilson-lines as representations
of groupoids (up to homotopy).

Let us note on our convention: the ground field in this paper, denoted as $k$, is either $\BB{R}$ or $\BB{C}$ throughout.
We choose to work in the $C^{\infty}$-setting, that is, all manifolds are taken to be smooth; even though we believe most
of the results can be transferred to the analytic and algebraic setting with minor modifications.
We work in the category of differential graded manifold where the parity is compatible with the grading
($\mathbb{Z}_2$- and  $\mathbb{Z}$-gradings), see \cite{voronov, Roytenberg:2002nu} for further
explanations and definitions.

\section{Representation of $NQ$-manifolds}\label{RoQS}

The data of a Lie algebroid can be easily given in terms of an $NQ$-manifold $\M$. Here $N$ refers to $\M$ being a non-negatively graded supermanifold, and the $Q$-structure refers to
a degree 1 nilpotent vector field, called the homological vector field or simply $Q$-vector field.
The $N$ and $Q$ structure are encoded by saying that there is an action of the semi-group
\bea \Hom(\BB R^{0|1},\BB R^{0|1})=\BB{R}^{1|1}\nn\eea
on $\M$ \cite{Severa06}. To see this, note that $\Hom(\BB R^{0|1},\BB R^{0|1})$ acts on $\BB R^{0|1}$ by shifting $\theta\to \theta+\epsilon$ and scaling $\theta\to s\theta,~s\in\BB{R}$. The two actions are generated by the vector fields
\bea \frac{\partial}{\partial\theta},~~~~\theta\frac{\partial}{\partial\theta}\nn\eea
The former vector field is clearly nilpotent and its action on $\M$ is the $Q$-vector field. The second gives the non-negative grading on $\M$.

\begin{Exa}\label{Lie_algebroid}
Let $M$ be a smooth manifold and $L\to M$ be a vector bundle, from which one can form the NQ-manifold $\M=L[1]$. We denote by $x^{\mu}$ the coordinate of $M$ and
$\ell^A$ the coordinate of the fibre. The grading of $x^{\mu},\ell^A$ is 0 and 1, and the functions $C^{\infty}(L[1])$ are isomorphic to  sections of  $\wedge^{\sbullet}L^*$.

The $Q$-structure, which is a deg 1 vector field over $\M$ is necessarily of the form
\bea Q=2\ell^AA_A^{\mu}(x) \frac{\partial}{\partial
x^{\mu}}-f^A_{BC}(x)
\ell^B\ell^C\frac{\partial}{\partial\ell^A}~,\label{Q_Lie_algebroid}\eea
where  $Q$ is written in local coordinates, requiring it to be globally defined will tell us how $A_A^{\mu}$ and $f_{BC}^A$ transform.
Imposing $Q^2=0$ puts constraint on the coefficients $A_A^{\mu}(x),f_{BC}^A(x)$
\bea&&A_{[A}^{\nu}\partial_{\nu}A_{B]}^{\mu}=A^{\mu}_Cf^C_{AB}~,\nn\\
&&A_A^{\mu}\partial_{\mu}f^D_{BC}+f^D_{AX}f^X_{BC}+\textrm{cyclic
in \scriptsize$ABC$}=0~,\label{Lie_Alge_compatible}\eea
$A_A^{\mu}$ is the called the 'anchor' giving a bundle morphism $\rho:\,L\to TM$
\bea \rho(s_{A})=A_A^{\mu}\frac{\partial}{\partial x_{\mu}},\nn\eea
where $s_A$ is a local basis of sections of $L$. And $f^A_{BC}$ is the 'structure
function' giving rise to a bracket on $\Gamma(L)$ as
\bea [s_A,s_B]=f_{AB}^Cs_C.\nn\eea
Then the first condition of Eq.\ref{Lie_Alge_compatible} says that the anchor $\rho$ respects the bracket,
\bea \rho([s_A,s_B])=[\rho(s_A),\rho(s_B)],\nn\eea
where the second bracket is the Lie bracket on $TM$. Whereas the second condition describes the failure of the Jacobi identity
\bea \rho(s_A)[s_B,s_C]+[s_A,[s_B,s_C]]+\textrm{cyclic
in \scriptsize$ABC$}=0\nn\eea
as measured by the anchor. These two conditions constitute the definition of a Lie algebroid; this formulation of the Lie algebroids is given in ref.\cite{vaintrob}.

The functions $C^{\infty}(L[1])$ are isomorphic to $\Gamma(\wedge^{\sbullet}L^*)$ which has the structure of the
 Chevalley-Eilenberg complex, and the $Q$-structure induces the Chevalley-Eilenberg differential; more of this will appear in sec.\ref{MCSaI}.

The shifted tangent bundle $T[1]M$ is the simplest example of a Lie algebroid, for this reason, Lie algebroids are sometimes called generalized tangent bundles.
\end{Exa}

We turn to the representations of Lie algebroids and also in more generality the representations of $NQ$-manifolds. Following Vaintrob \cite{vaintrob}, for an $NQ$-manifold $\M$, consider a graded vector bundle\footnote{Graded vector bundles over graded manifolds are defined as sheaves of freely generated $C^{\infty}(\M)$-modules over the reduced manifold of $\M$, the coordinates of the fibre in the following discussion are the generators of this module, see ref.\cite{NoteonSusy}.} $\E$ over ${\cal M}$.
\bea
\begin{tikzpicture}[scale=.8]\label{graded_bundle}
  \matrix (m) [matrix of math nodes, row sep=2.5em, column sep=2.5em]
    { \hat Q & {\cal E} & {\cal F}\\
      Q & {\cal M} &  \\ };
 \path[->]
 (m-1-3) edge  (m-1-2)
 (m-1-2) edge  node[right] {$\small{\pi}$} (m-2-2)
 (m-1-1) edge  node[right] {$\small{\pi}_*$} (m-2-1);
\end{tikzpicture}\label{def_rep}
\eea
Further let us assume that $\E$ is $Q$-equivariant, meaning the vector field $Q$ on $\M$ can be lifted to $\hat Q$ over $\E$ with $\hat{Q}^2=0$. In other words there exists a $Q$-flat connection $\R\in\textrm{End}_{C^{\infty}(\M)}(\E)$ such that
\bea \hat Q^2=(Q+\R)^2=0.\nn\eea
The $Q$-flat connection $\R$ is the \emph{representation of the $Q$-structure}. We shall see in later examples that this definition encodes the representation of Lie algebroids up to homotopy very compactly.

We can choose a trivialization locally on a patch $\E|_U=\M|_U\times {\cal F}$, and denote the coordinate of the fibre ${\cal F}$ as $\zeta_{\alpha}$ and that of $\M$ as $x^A$. Locally, we can write $\hat{Q}$ as
\bea \hat Q=Q^A(x)\frac{\partial}{\partial x^A}+R^{\beta}_{~\alpha}(x)\zeta_{\beta}\frac{\partial}{\partial \zeta_{\alpha}},\label{representation}\eea
where the fibre dependence of $\R$ is by definition linear. The nilpotency of $\hat Q$ says
\bea QR^{\beta}_{~\alpha}+(-1)^{(\beta+\gamma)(1+\gamma+\alpha)}R^{\beta}_{~\gamma}R^{\gamma}_{~\alpha}=0,\nn\eea
where the degree of $\zeta_{\ga}$ is denoted as $\ga$ and the degree of $R^{\gamma}_{~\alpha}$ is $1+\gamma+\alpha$. Let us define the normalized representation matrix
\bea&& T^{\beta}_{~\alpha}=(-1)^{\alpha\beta+\beta}R^{\beta}_{~\alpha}~~\Rightarrow~~ (-1)^{\beta}QT^{\beta}_{~\alpha}+(-1)^{\gamma}T^{\beta}_{~\gamma}T^{\gamma}_{~\alpha}=0\label{normalized_MC},\eea
it will be clear in a moment what is the point of such redefinition.

As a matter of notation we will call the degree of the coordinates of $\M$ \emph{degree} and the grading of the fibre coordinates of $\E$ \emph{level}. Thus the matrix $T^{\ga}_{~\gb}$ shifts the level up by $\ga-\gb$. For a non-negatively graded manifold, we have the decomposition of $T$ according to its degree
\bea T^{\ga}_{~\gb}=T^{\ga}_{0\gb}+T^{\ga}_{1\gb}+\cdots\label{decomp_T},~~~~\deg T^{\ga}_{p\gb}=p\label{decomp_T}.\eea
The following is a trivial consequence of the decomposition
\bea T^{\ga}_{p\gb}=\Big\{\begin{array}{cc}
                       T^{\ga}_{~\gb} & \ga=\gb-p+1\\
                       0 & \textrm{otherwise}
                     \end{array}\nn.\eea
In particular, the term $T_{0}$ raises the level of $\E$ by 1, and we have as a consequence of Eq.\ref{normalized_MC}
\bea T_0^2=0.\nn \label{T-square}
\eea
For this very reason, $T_0$ is a differential that moves one up in the fibre level, leading to a fibrewise structure of a differential chain complex. Later we shall see in the phrase 'representation up to homotopy', the homotopy is the chain homotopy of $T_0$ (or in other words, 'things work just fine' up to $T_0$-exact terms).

\begin{Exa}\label{EX_Flat_bundle} Flat bundle as a representation of $T[1]M$

Take $\M=T[1]M$, where $TM$ is the tangent bundle of a smooth manifold $M$, and the coordinate of the fibre is assigned degree 1. So if we let $x^{\mu}$ be the local coordinates of $M$ and $v^{\mu}$ be the deg 1 coordinate of the fibre, then the (sheaf of) functions on $T[1]M$ are
\bea C^{\infty}(T[1]M)= \Omega^{\sbullet}(M),\nn\eea
in other words $v^{\mu}$ are identified as $dx^{\mu}$.

The Lie algebroid data is given by the $Q$-structure on the $T[1]M$ which corresponds to the de Rham differential
\bea D=v^{\mu}\frac{\partial}{\partial x^{\mu}}.\nn\eea

The graded vector bundle over $T[1]M$ is taken to be the pull back of a flat bundle
$E\to M$ to $T[1]M$. Denoting the flat connection of $E$ as $A$, then
we can write down a lift of $D$ as
\bea \hat{D}=v^{\mu}\frac{\partial}{\partial x^{\mu}}+v^{\mu}(A_{\mu})^{i}_{~j}\zeta^{j}\frac{\partial}{\partial \zeta^{i}}.\label{Q_flat_bundle}\eea
The nilpotency of $\hat D$ is equivalent to the flatness of $A$.

Our construction seems to have excluded non-flat vector bundles. But in fact even if $E$ is non-flat one can write down a lift of $D$ by going to the jet bundle $E\otimes \textrm{jet}^{\infty}(TM)$, however we shall merely give the formula and not digress too far afield
\bea \hat D=v^{\ga}\partial_{\ga}+v^{\gc}\big(\delta^{\mu}_{\gc}+\Gamma^{\mu}_{\gb\gc}\xi^{\gb}-\frac13R_{\gc\ga~\gb}^{~~\;\mu}\xi^{\ga}\xi^{\gb}+\cdots\big)\frac{\partial}{\partial \xi^{\mu}}+v^{\gc}\big(A_{\gc}+\frac12\xi^{\gb}F_{\gb\gc}+\frac16\xi^{\ga}\xi^{\gb}\BS{\nabla}_{\ga}F_{\gc\gb}+\cdots\big)^i_{~k}\zeta^k\frac{\partial}{\partial \zeta^{i}},\nn\eea
where $\xi^{\ga}$ is the coordinate of $TM$ (so the section of $\textrm{jet}^{\infty}(TM)$ is a formal polynomial of $\xi$), and $F$, $R$ are the curvature of $E$ and $TM$. By restricting to the subbundle $E$, we get the Wilson-line for the non-flat connection $A$.
\end{Exa}

\begin{Exa}\label{Courant} Standard Courant Algebroid as the adjoint representation of $T[1]M$

The Lie algebroid data is given as Ex.\ref{EX_Flat_bundle}, and the graded vector bundle over $T[1]M$ is taken to be
\bea \E=T^*[2]T[1]M,\nn\eea
and locally we denote the fibre coordinate of $\E$ as $p_{\mu},q_{\mu}$ of degree 2,1 and they are the 'momentum' dual to $x^{\mu}$, $v^{\mu}$ respectively. We can write down a lift of $D$ as
\bea \hat{D}=v^{\mu}\frac{\partial}{\partial x^{\mu}}+p_{\mu}\frac{\partial}{\partial q_{\mu}}.\label{D_hat}\eea
Upon choosing a connection and split $T^*[2]T[1]M$ into $(T^*[2]\oplus T^*[1]\oplus T[1])M$, we recover the adjoint representation up to homotopy of $T[1]M$ in the sense of ref.\cite{abad-2009}. In this example, only $T_0$ is non-zero. The construction using a Hamiltonian lift is completely general, we provide more details in the next example.
\end{Exa}

\begin{Exa}\label{Adjoint_rep} Adjoint representation (up to homotopy) of a Lie algebroid

The Lie algebroid data is given by the $Q$-structure as in Ex.\ref{Lie_algebroid}. And consider the graded vector bundle over $L[1]$
\bea \E=T^*[2]L[1],\nn\eea
and locally we denote the fibre coordinate of $\E$ as $p_{\mu},\bar\ell_A$ of degree 2,1, which are the momentum dual to $x^{\mu}$, $\ell^A$ respectively. We can write down a lift of $Q$ using the following recipe.
Since $T^*[2]L[1]$ is symplectic with the standard symplectic form
\bea \Omega=dp_{\mu}dx^{\mu}+d\bar\ell_Ad\ell^A\nn\eea
of degree 2. Let us lift $Q$ up into a Hamiltonian vector field by replacing $\partial_{\mu}$ with $p_{\mu}$ and $\partial_A$ with $\bar\ell_A$
\bea Q\to \Theta=2p_{\mu}A^{\mu}_A\ell^A+f^A_{BC}\bar\ell_A\ell^B\ell^C\nn\eea
The lift of $Q$ is defined as the Hamiltonian vector field generated by $\Theta$
\bea \hat{Q}\cdot=\{\Theta,\cdot\}=Q+\cdots.\nn\eea
Clearly the terms $\cdots$ are linear in the fibre coordinates $p_{\mu}$ and $\bar\ell_A$ as $\textrm{End}_{C^{\infty}(\M)}\E$ should be. We show in the appendix how upon picking a connection for $L$ one recovers the notion of representation of Lie algebroids up to homotopy.
\end{Exa}

\section{Superfields and the Functor of Points Perspective}\label{sec_SF}
We review briefly the 'functor of points' view on the graded manifolds and explain how this is just the superfield technique that the physicists have been using all along.

The graded manifolds (GM) are locally ringed spaces, that is, a graded manifold $\M$ consists of a smooth manifold $M$ with a sheaf of freely generated graded commutative algebra on $M$. The \emph{section}s of this sheaf, denoted $C^{\infty}(\M)$, is often called the \emph{function}s on the graded manifold by a standard abuse of terminology (see ref.\cite{NoteonSusy}). The underlying smooth manifold $M$ is called the reduced manifold of $\M$ denoted as $M=|\M|$ ($M$ is also called the body of $\M$). Morphisms (or simply, maps) between graded manifolds are morphisms of their sheafs of algebras, $\Hom_{GM}(\M,\N)\sim\Hom_{\gam}(C^{\infty}(\N), C^{\infty}(\M))$, where $\Hom_{\gam}$ stands for morphisms of graded commutative algebras. We will use the term 'GM maps' and '\gam-maps' interchangeably.

The ringed structure makes it awkward to talk about points on a GM, a more efficient way of dealing with GM's is to look at the structure of the algebra of functions, and the morphisms of functions. For the category of GM's, one can define a presheaf which is a contravariant functor from GM to the category of sets. The GM's themselves give such a functor: if $\M$ is a GM it gives a presheaf $\hat{\M}$ defined as \bea \hat{\M}:~\textrm{GM}\longrightarrow \Set,~~~~~\hat{\M}(\Pa)=\Hom_{GM}(\Pa,\M),\label{presheaf_hat}\eea where $\Pa$ is a GM and $\Hom_{GM}(\Pa,\M)$ is the set of morphisms from $\Pa$ to $\M$. The presheafs of the above type are called \emph{representable}. Identifying $\M$ as such a presheaf is the 'functor of points' way of perceiving a GM: for each $\Pa$, $\hat\M(\Pa)$ is a family of points parameterized by $\Pa$, in other words we use a 'test manifold' $\Pa$ to probe $\M$. For example the set $\Hom_{GM}(*,\M)$ gives the reduced manifold of $\M$. And Yoneda's lemma guarantees that the set of natural transformations between $\hat{\M}$ and $\hat{\N}$ is $\Hom_{GM}(\M,\N)$. This point of view also allows us to think about slightly more general GM's as presheafs that are not necessarily representable.

Another important example of a presheaf is
\bea \HomM_{GM}(\So,\Ta): \textrm{GM}\longrightarrow \Set,~~~~~\HomM_{GM}(\So,\Ta)({\Pa})=\Hom_{GM}(\So\times\Pa,\Ta),\label{presheaf_hom}\eea
where $\Pa$ is a GM and $\Hom(\So\times\Pa,\Ta)$ is the set of morphisms $\So\times\Pa\to\Ta$. Of course, both presheafs Eq.\ref{presheaf_hat}, \ref{presheaf_hom} are trivially generalizable to the category of $NQ$-manifolds.

It is known that the presheaf Eq.\ref{presheaf_hom} is representable, and the GM that represents it is called \emph{superfields}.
Next we introduce the superfields as an important computational tool and the discussion here largely follows ref.\cite{Roytenberg:2006qz}.
The statement is near vacuous, to see this, let us denote the algebra of functions on $\Ta$, $\So$, $\Pa$ as $\textsf{A}$, $\textsf B$ and $\textsf{C}$, then a morphism $\So\to\Ta$ (resp. $\So\times\Pa\to\Ta)$ is \gam-map $\textsf{A}$ to $\textsf{B}$ (resp. $\textsf{A}\to\textsf{B}\otimes\textsf{C}$). Let
$x^m$, $y^a$ and $z^i$ be the generators of $\textsf{A, B}$ and $\textsf{C}$ respectively, and we also use a subscript to denote the degree of the generators. Now one can expand a map
\bea \textsf{A}\stackrel{\varphi}{\longrightarrow}\textsf{B}\otimes\textsf{C}:~x^m\to \BS{x}^m(y,z),\nn\eea
as a power series of $y$ with non-zero degrees (for clarity, we assume $\So$ is non-negatively graded)
\bea &&\BS{x}_p(y,z)=x_{p,p}(y_0,z)+x_{p,p-1}(y_0,z)y_1+\Big(\frac{1}{2}x_{p,p-2}(y_0,z)y_1^2+x^{\prime}_{p,p-2}(y_0,z)y_2\Big)\nn\\
&&\hspace{2cm}+\Big(\frac16x_{p,p-3}(y_0,z)y_1^3+x^{\prime}_{p,p-3}(y_0,z)y_1y_2+x^{\prime\prime}_{p,p-3}(y_0,z)y_3\Big)+\cdots,\label{SF_expansion}\eea
where we have suppressed the index structure for $x,~y$, but the meaning should be clear. The notation $x_{p,q}$ denotes the degree $q$ component of $\BS{x}_p(y,z)$ in the expansion. Notice that we have only expanded $\BS{x}(y,z)$ w.r.t the $y$'s and kept the $z$'s intact. One can take $x^m_{p,q}$ as the generators of the algebra of functions of an infinite dimensional graded manifold which we call ${\cal S}$-the superfields, then the set of maps $\Hom_{GM}(\So\times\Pa,\Ta)$ is given by the set of maps $\Hom_{GM}(\Pa,{\cal S})$, leading to the statement about representability. The arguments can be repeated for the case of $N$-manifolds except that the expansion in Eq.\ref{SF_expansion} will terminate at $x^{m}_{p,0}$.

If we set the test manifold $\Pa$ to be a point, then $\textsf{C}=k$, as a result only the degree zero component $x_{p,0}$ in the expansion survives\footnote{This is simply because there are no $\deg>0$ generators left in town to write down a component $x_{p,q},\;q>0$}
\bea \BS{x}_p(y)=\frac1{p!}x_{p,0}(y_0)y_1^p+\frac1{(p-2)!}x^{\prime}_{p,0}(y_0)y^{p-2}_1y_2+\cdots+x^{\prime\cdots\prime}_{p,0}(y_0)y_p.\nn\eea
Except in our discussion of homotopies of \cdga-maps $\Hom_{\cdga}(\textsf{A,B})$ coming up shortly, we can, for most calculations, be utterly oblivious to what exactly is $\L$ or what constitutes $\textsf{C}$, and write the expansion Eq.\ref{SF_expansion} simply as
\bea \BS{x}_p(y)=x_{p,p}(y_0)+x_{p,p-1}(y_0)y_1+\frac12x_{p,p-2}(y_0)y_1^2+x^{\prime}_{p,p-2}(y_0)y_2+\cdots\nn,\eea
keeping all the components $x_{p,q}(y_0)$. Since we will only be interested in $N$-manifolds, only $x_{p,\geq0}(y_0)$ should be kept as {\sf C} is now non-negatively graded, but keeping the negative degree components lubricates the calculations appreciably, and as such, the negative degree components are called the 'auxiliary fields', in keeping with their name in physics. The negative degree components will be set to zero as soon as the calculation is done.

Up to now, one may well wonder what on earth is novel about this point of view. Indeed, the above discussion shows we can carry most of our intuition in differential geometry to the graded case and treat the generators $x^m$ $y^a$ as if they were coordinates and can take specific values. The physicists have carried out many brilliant computations in such way without ever bothering about the category jargons.

\subsection{$\cdga$-Maps}
So far we have discussed graded algebra maps, now if $\textsf{A,B}$ and $\textsf{C}$ are actually (non-negatively) graded commutative differential algebras (\cdga's) equipped with the differentials $\gd,d^y$ and $d^z$ respectively, 
then the condition for any \gam-map $\varphi:~ x^m\to \BS{x}^m(y,z)$ to be a \cdga-map is
\bea 0=(\gd \BS{x}^m)(y,z)-(d^y+d^z)(\BS{x}^m(y,z)),\label{condition_dga}\eea
namely $\varphi$ commutes with the differentials. By looking at the equation at $z_{>0}$ all equal to zero, we have
\bea 0=(\gd \BS{x}^m)(y,z_0)-d^y(\BS{x}^m(y,z_0)),\nn\eea
which means $\BS{x}^m(y,z_0)$ is a \cdga-map $\textsf{A}\to\textsf{B}$ for all $z_0$. Recall that once we set $z_{>0}=0$, $x_{p,>0}$ are also zero.
We denote (suppressing $z$)
\bea \textrm{eom}=(\gd \BS{x}^m)(y)-d^y(\BS{x}^m(y))\Big|_{\deg\,0}.\nn\eea
as the \emph{equation of motion}.

There is a $Q$-structure defined on ${\cal S}$, the \emph{BRST differential}, acting on the components $x^m_{p,q}(y_0)$, induced by $\gd$ and $d^y$
\bea (\delta_B\BS{x}^m)(y)&=&\gd_Bx^m_{p,p}(y_0)+(\gd_Bx^m_{p,p-1}(y_0))y_1+\frac{1}{2}(\gd_Bx^m_{p,p-2}(y_0))\frac12y^2_1+(\gd_Bx^{\prime m}_{p,p-2}(y_0))y_2+\cdots\nn\\
&=&(\gd \BS{x}^m)(y)-d^y(\BS{x}^m(y)).\label{BRST_abstract}\eea
By comparing the power of $y$ on both sides, one can read off the action of $\gd_B$ on each component\footnote{We would like to stress that $\delta_B$ acts on the components fields $x_{p,q}(y_0)$, not on $x$ nor $y$, and only when the fields are assembled into superfields, can the action be written nicely as
$(\delta_B\BS{x}^m)(y)=(\gd \BS{x}^m)(y)-d^y(\BS{x}^m(y))$.}.
The differential $\gd_B$ plays the crucial role in this paper. The degree zero component of Eq.\ref{BRST_abstract} is the eom defined earlier and quite clearly, the equation of motion is the BRST differential of the deg $-1$ component
\bea \textrm{eom for }x^m=\delta_B x^m_{p,-1}\Big|_{\deg\;0}.\label{aux}\eea

We introduce a space $\bar{\cal S}\subset{\cal S}$, whose algebra of functions is
\bea C^{\infty}(\bar{\cal S})=C^{\infty}({\cal S})\big/(\textrm{eom},x^m_{p,<0}(y_0)),\label{ideal}\eea
where $(\textrm{eom},x^m_{p,<0}(y_0))$ is the ideal generated by eom and $x_{p,<0}$. Property Eq.\ref{aux}
shows that the quotient is compatible with $\gd_B$.

And from Eq.\ref{condition_dga} the following is clear
\begin{Prop}\label{representability}
$\bar{\cal S}$ is the $NQ$-manifold that represents the presheaf $\HomM_{NQ}(\So,\Ta)$
\end{Prop}
We remind the reader that $\HomM_{NQ}$ is a presheaf defined similarly to Eq.\ref{presheaf_hom}, but with $GM$ replaced with $NQ$.

The compatibility of $\gd_B$ with the quotient also allows us do all the computation in the space ${\cal S}$, and only in the end go over to $\bar{\cal S}$.

\subsection{Homotopy of \cdga-Map}
We now follow ref.\cite{DGMS} (see also ref.\cite{Sull}) in the formulation of homotopy of \cdga-maps. To deform a \cdga-map $\A\to\textsf{B}$, we take $\textsf{C}$ as $C^{\infty}(T[1][0,1])$ and thus $\textsf{C}=k[s,\vartheta]$ with $\vartheta$ being the odd coordinate on $T[1][0,1]$. The differential $d^s$ is defined in the obvious way $d^ss=\vartheta$, $d^s\vartheta=0$. Given two \cdga-maps $\textsf{A}\stackrel{\varphi,\phi}{\to}\textsf B$, the homotopy is a \cdga-map
\bea \textsf{A}\stackrel{H}\to \textsf{B}[s,\vartheta],\label{homotopy_setting}\eea
such that $H$ at $s=0,1$ equals $\varphi$, $\phi$ respectively. In terms of coordinates, any map $H$ is written as $x^m\to \BS{x}^m(y,s,\vartheta)$ and the condition for \cdga-map is
\bea0=(\gd \BS{x}^m)(y,s,\vartheta)-(d^y+d^s)(\BS{x}^m(y,s,\vartheta)).\label{dga_homotopy}\eea
If we partially expand the superfield as
\bea \BS{x}^m(y,s,\vartheta)=\BS{x}^m(y,s)+\vartheta\bar{\BS{x}}^m(y,s),\nn\eea
where the expansion terminates because $\vartheta^2=0$ trivially and we also have $\deg \bar x^m=\deg x^m-1$.
Then looking at Eq.\ref{dga_homotopy} at $\vartheta=0$ we get
\bea 0=(\gd\BS{x}^m)(y,s)-d^y(\BS{x}^m(y,s)),\label{eom_abstract}\eea
i.e. $\BS{x}^m(y,s)$ is a \cdga-map $\A\to\textsf{B}$ for each $s$. Looking at the same equation at first power of $\vartheta$ we get
\bea 0=\vartheta\Big((\bar{\BS{x}}^n\partial_{x^n}\gd\BS{x}^m)(y,s)-\partial_s \BS{x}^m(y,s)+d^y(\bar{\BS{x}}^m(y,s))\Big),\label{name_later}\eea
giving the formula for an infinitesimal homotopy
\bea \partial_s \BS{x}^m(y,s)=(\bar{\BS{x}}^n\partial_{x^n}\gd\BS{x}^m)(y,s)+d^y(\bar{\BS{x}}^m(y,s)).\label{homotopy_abstract}\eea

%

\begin{rmk}
We can write a homotopy as if it were a Lie derivative given by the Cartan formula
\bea \frac{\partial}{\partial s}\BS{x}^m(y,s)=\L_{\bar x} \BS{x}^m(y)=\big\{\iota_{\bar x}, \delta_B\big\}\BS{x}^m(y),~~~~~\iota_{\bar x}=\bar x^m \frac{\partial}{\partial x^m}.\label{homotopy_cartan}\eea
\end{rmk}
The homotopies Eq.\ref{homotopy_abstract} are called \emph{gauge transformations} in physics.
We give some examples next.

\subsection{A Lie Algebroid Illustration}\label{ex_Lie alge_eom}
We look at the \cdga-map $\Hom_{\cdga}\big(C^{\infty}(L[1]),C^{\infty}(T[1][0,1])\big)=\Hom_{NQ}\big(T[1][0,1],L[1]\big)$.

Referring back to Ex.\ref{Lie_algebroid}, we let $\textsf{A}=(C^{\infty}(L[1]),Q)$. We also let $\textsf{B}=(C^{\infty}(T[1][0,1]),D)$, concretely we have $\textsf{B}=k[t,\theta]$, with $D=\theta\partial_t$.
We have the superfield expansion (keeping all components for the time being, namely we work in ${\cal S}$ for now)
\bea \BS{x}^{\mu}(t,\theta)=x^{\mu}(t)+\theta\,x^{\mu}_t(t),~~~~~~\BS{\ell}^A(t,\theta)=\ell^A(t)+\theta\,\ell^A_t(t).\nn\eea
Note $x^{\mu}_t$ is of deg $-1$.

We can read off $\gd_B$ on the various components from Eq.\ref{BRST_abstract}
\bea (\gd_B\BS{x}^{\mu})=2A^{\mu}_A(\BS{x})\BS{\ell}^A-D\BS{x}^{\mu},
~~~~(\gd_B\BS{\ell}^A)=-f^A_{BC}(\BS{x})\BS{\ell}^B\BS{\ell}^AD\BS{\ell}^A\nn\eea
as
\bea \gd_Bx^{\mu}(t)=2A^{\mu}_A(x(t))\ell^A(t),&&\gd_Bx^{\mu}_t(t)=-2x^{\nu}_t\big(\partial_{\nu}A^{\mu}_A(x(t))\big)\ell^A(t)-2A^{\mu}_A(x(t))\ell^A_t(t)+\dot x^{\mu}(t),\nn\\
\delta_B\ell^A(t)=-f^A_{BC}\ell^B(t)\ell^C(t),&&\delta_B\ell^A_t(t)=x^{\mu}_t(t)\big(\partial_{\mu}f^A_{BC}(x(t))\big)\ell^B(t)\ell^C(t)
+2f^A_{BC}(x(t))\ell_t^B(t)\ell^C(t)+\dot\ell^A(t),\nn\eea
where a dot denotes the derivative with respect to $t$.

Dropping the cumbersome argument $t$ everywhere, the equation of motion is
\bea 0=\gd_Bx^{\mu}_t\Big|_{\deg\,0}=\partial_tx^{\mu}-2A^{\mu}_A(x)\ell^A_t,\label{eom_lie_algebroid}\eea
which actually is the condition that the mapping (where $e^A$ is the basis of the sections of $CE^1$)
\bea x^{\mu}\to x^{\mu}(t),~~~e^A\to \ell^A_t(t) \theta,\nn\eea
is a morphism of \cdga~ $CE^{\sbullet}\to \Omega^{\sbullet}([0,1])$.

Specializing to $L[1]=T[1]M$ (see Ex.\ref{EX_Flat_bundle} for the notation), the above equation of motion reads
\bea \dot x^{\mu}=v^{\mu}_t,\nn\eea
hence we can eliminate $v_t^{\mu}$ in favour of $\partial_tx^{\mu}$. Thus we get the trivial statement: any \cdga-map from the de Rham complex of $M$ to that of the interval is completely determined by a smooth map $[0,1]\to M$ and the set $\Hom_{NQ}(T[1][0,1],T[1]M)$ is the path space $PM$.

Likewise, a homotopy is fixed by a tuple $(\bar x^{\mu},\bar\ell^A)$. For clarity, take the test manifold to be a point, then the abstract formula Eq.\ref{homotopy_abstract} unfolds as
\bea &&\partial_sx^{\mu}=2A^{\mu}_A\bar\ell^A,\nn\\
&&  \partial_s\ell^A_t =-2f^A_{BC}\bar\ell^B\ell_t^C+\dot{\bar\ell}^A.\label{homotopy_liealbegroid}\eea
In the calculation we have taken the short cut by only keeping the degree 0 components, for example the first equation should really go as
\bea &&\partial_s x^{\mu}=\bar x^{\nu}(\partial_{\nu}A_A^{\mu})\ell^A+2A^{\mu}_A\bar\ell^A,\nn\eea
yet only the second term survives when the test manifold is a point.

The homotopy Eq.\ref{homotopy_liealbegroid} by construction takes a \cdga-map to a \cdga-map.
One can of course check this explicitly
\bea 2\partial_s\big(\ell^A_tA_A^{\mu}\big)&=&2\big(\dot{\bar\ell}^A-2{\bar\ell}^Bf_{BC}^A\ell^C_t\big)A_A^{\mu}+\ell^A_t\big(2{\bar\ell}^BA_B^{\nu}\partial_{\nu}A_A^{\mu}\big)\nn\\
&=&2\dot{\bar\ell}^AA_A^{\mu}+4\ell^A_t{\bar\ell}^BA_A^{\nu}\partial_{\nu}A_B^{\mu}=2\dot{\bar\ell}^AA_A^{\mu}+2{\bar\ell}^B\dot x^{\nu}\partial_{\nu}A_B^{\mu}=\partial_s\dot x^{\mu}\nn,\eea
which shows that Eq.\ref{eom_lie_algebroid} is satisfied for every $s$.

Before leaving this example, we point out that for the case $L[1]=T[1]M$, the formula Eq.\ref{homotopy_cartan} is nothing but the Cartan formula for Lie derivatives.


 We would like to point out the use of homotopy in this paper refers to homotopies of \cdga-maps; in the case of \cdga~maps between $C^{\infty}(T[1]M)$ and $C^{\infty}(T[1]N)$, then the homotopy of the \cdga~map is indeed the homotopy of maps $N\to M$ in the conventional sense. But in general, the homotopy can be a bit unconventional, in fact two \cdga-maps that are 'near each other' may not always be homotopic.
\begin{Exa}Consider the space $\BB{R}^3$ with the Poisson tensor
\bea \ga^{\mu\nu}=x^{\rho}\epsilon^{\mu\nu\rho},\nn\eea
and let $L=T^*\BB{R}^3$. One can form the Lie algebroid by specifying
\bea Q=2\alpha^{\mu\nu}\xi_{\mu}\frac{\partial}{\partial
x^{\nu}}+(\frac{\partial}{\partial
x^{\rho}}\alpha^{\mu\nu})\xi_{\mu}\xi_{\nu}\frac{\partial}{\partial
\xi_{\rho}},\nn\eea
where $x,\xi$ are the degree 0,1 coordinates of $T^*[1]\BB{R}^3$. The condition for \cdga-map implies
\bea \dot x^{\mu}(t)=\ga^{\mu\nu}(x(t))\xi_{t\nu}(t).\nn\eea
This condition says that the image of the interval $[0,1]$ is contained in the leaf of the Poisson tensor, which for our case are just concentric spheres.

If we look at the homotopy, these are deformations of the form Eq.\ref{homotopy_liealbegroid}
\bea
\partial_sx^{\mu}=\bar\xi_{\nu}\ga^{\mu\nu},\eea
which again must land in the same leaf. Thus two curves $[0,1]\to\BB{R}^3$ belonging to two different leafs are not homotopic in our sense, no matter how close.\end{Exa}

\section{Maurer-Cartan Set and Integration}\label{MCSaI}
We first define the Chevalley-Eilenberg complex of a Lie algebroid $L\to M$. Let $s_1,s_2,\cdots$ be sections of $L$, and $c^q$ is an antisymmetric linear functional on sections of $L$
\bea &&c^q\in CE^q,~c^q(s_1,\cdots s_q)\in C^{\infty}(M),\nn\eea
The complex has a differential
\bea \delta c^q(s_1,\cdots s_{q+1})=\sum_{i=1}^q \rho(s_i)\circ c^q(s_1,\cdots \hat{s_i},\cdots s_{q+1})
-\sum_{i<j}(-1)^{i+j-1}c^q([s_i,s_j],s_1,\cdots \hat{s_i},\cdots\hat{s_j},\cdots s_{q+1}),\label{CE_def}\eea
where $\rho(s)\circ f(x)=s^AA_A^{\mu}\partial_{\mu}f(x)$. Just like the case of differential forms, $c^q$ is $C^{\infty}(M)$-linear in all of its entries.

The linearity implies that the Chevalley-Eilenberg complex is spanned by sections of the graded vector bundle
\bea \E=\oplus_q\wedge^q L^*,\nn\eea
and as a general rule in this paper, \emph{the basis of the sections of a graded bundle is regarded as the coordinates of the dual bundle with its degree equal to the grading of the section in question}.
Thus the coordinates $\ell^A$ of $L[1]$ represents the basis of sections of $L^*$ regarded as the degree 1 element of $CE^1$, and $CE^q$ is spanned by the functions $C^{\infty}(L[1])$ of degree $q$.
The CE differential Eq.\ref{CE_def} is simply given by the $Q$-vector field of Eq.\ref{Q_Lie_algebroid}.
This is probably the fastest way to see why the differential Eq.\ref{CE_def} is nilpotent in the first place. Taking $L[1]=\FR{g}[1]$, then $\rho=0$ and we recover the familiar CE complex for a Lie algebra.

Next let us review quickly the integration problem. Let $G$ be a Lie group whose Lie algebra is $\FR{g}$, then
locally, the '$G$-descent data', which is a map $M\times M\stackrel{\varphi}{\to} G$, such that
\bea \varphi(p,p)=e,~~~~\varphi(p,q)\varphi(q,s)=\varphi(p,s),~~~p,q,s\in U\subset M\label{G_descent}\eea
can be specified by a flat connection
\bea A\in \Omega^1(M)\otimes \FR{g},~~~~dA+A\wedge A=0.\nn\eea
To construct the map $\varphi$, for points $p,q\in U$, one picks a path $C$ from $p$ to $q$, and then one writes down the parallel transport (Wilson-line) operator
as a path ordered exponential
\bea (p,q)\rightarrow\BB{P}\exp\big\{\int_C\,A\big\},\nn\eea
where $\BB{P}$ denotes path ordering. The exponential can be identified as an element in the group $G$ by means of the Baker-Campbell-Hausdorf formula\footnote{Or one can pick a faithful representation for $\FR{g}$ and think of the exponential as a matrix exponentiation.}. The flatness condition ensures that the choice of path $C$ is immaterial and we get a well defined map from $(p,q)$ to $G$ with all the desired properties. This formulation of the mapping from $M$ to $G$ of course suffers from global issues, and we will not take on this problem here, the reader may see ref.\cite{Olver} for some cute examples of local integration that lacks global associativity and also consult ref.\cite{Zhu} for a proposal for the globalization procedure.

It is well known that a flat connection induces a chain map between $CE^{\sbullet}(\FR{g})$ and $\Omega^{\sbullet}(M)$, so the integration problem of a Lie algebra is to look for \cdga-maps $\Hom_{\cdga}\big(CE^{\sbullet}(\FR{g}),\Omega^{\sbullet}(M)\big)$.
As a straightforward generalization, for a Lie algebroid one looks for maps
\bea \Hom_{\cdga}\big(CE^{\sbullet}(L),\Omega^{\sbullet}(M)\big),\label{integration_Lie}\eea
which are called Maurer-Cartan elements.

Following our tenet in the italic font above, we recast Eq.\ref{integration_Lie} as an $NQ$ manifold map
\bea
\Hom_{NQ}\big(T[1]M,L[1]).\nn\eea
In fact, one can do the same for the integration problem of a general $NQ$-manifold. Let $\Ta$ be an $NQ$-manifold and $\So$ be a smooth manifold, we can define a presheaf $\int\M$, sending Man to the category of sets
\bea \smalint \Ta:~\textrm{Man}\to\Set,~~~~\smalint{\Ta}(\So)=\Hom_{NQ}(T[1]\So, \Ta),\label{significant}\eea
where $\Hom_{NQ}(T[1]\So, \Ta)$ is regarded as the set of $NQ$-morphism $T[1]\So\to\Ta$; this set is called the \emph{Maurer-Cartan set}. The reader may see ref.\cite{Severa06} for a nice account of how to view this functor as the adjoint of a 'differential functor' (called 1-jet in that paper).
Since by taking $\M=\FR{g}[1]$ and $\So=M$, the Maurer-Cartan set is the set of flat connections, it is natural to call the Maurer-Cartan set Eq,\ref{significant}
'the generalized flat connections'. We have seen that the formulation of integration in terms of CE complex is equivalent to the one in terms of the $NQ$-manifolds, in fact, from the latter point of view,
the integration is simply to integrate the flow of the homological vector field $Q$ on the $NQ$-manifold, which is very natural.

Note that in this paper, we aspire not to solve the integration problem, what we shall do is to construct certain functions on the presheaf $\Hom_{NQ}(T[1]\So, \Ta)$ that are invariant under the homotopies.

\section{Bar Complex and Iterated Integral}\label{BCaII}
We first give the definition of the bar resolution $B^{\cdot}(\A,\A,\Mod)$ of a differential graded module $\Mod$ over a \cdga~$\A$ with unit (for a more pedagogical review for connected $\A$ see ref.\cite{McCleary} ch.7).

The bar complex is a resolution of $\Mod$,
\bea \cdots \to B^{-q}\stackrel{b_1}{\to} B^{-q+1}\stackrel{b_1}{\to}\cdots B^0\stackrel{\epsilon}{\to} \Mod\to 0,\nn\eea
at degree $-q$ it is defined as
\bea B^{-q}(\A,\A,\Mod)=\A\otimes_k\underbrace{\A\otimes_k\A\cdots\A}_q\otimes_k \Mod,\label{def_bar_un}\eea
and a typical element is denoted as
\bea f_0[f_1|f_2|\cdots |f_q]m\in B^{-q}(\A,\A,\Mod),\nn\eea
with $f_i\in \A$ and $m\in\Mod$. The augmentation and the embedding is defined
\bea &&\epsilon:~ B^{0}\to \Mod,~~~\epsilon(f_0[~]m)=f_0m,\nn\\
&&\eta:~ \Mod \to B^{0},~~~\eta(m)=1[~]m.\nn\eea
And $B^{\sbullet}$ has a differential $b_1$ that acts on a typical element as
\bea b_1\Big(f_0\left[f_1|f_2|\cdots|f_q\right]m\Big)&=&(-1)^{f_0}f_0f_1[f_2|f_3|\cdots|f_q]m\nn\\
&&+\sum_{i=1}^{q-1}f_0[f_1|\cdots |f_{i-1}|f_if_{i+1}|f_{i+2}|\cdots|f_q]m(-1)^{f_0+(f_1+1)+\cdots +(f_i+1)}\nn\\
&&-f_0[f_1|\cdots|f_{q-1}]f_qm(-1)^{f_0+(f_1+1)+\cdots +(f_{q-1}+1)}.\label{diff_bar}\eea
where $(-1)^{f_0}$ means $(-1)^{\deg\,f_0}$ and so forth. Rather than checking now the property $b^2_1=0$, it will become clear after Eq.\ref{bar_inter_it}.
Besides, $b_1$ is acyclic, since one can define a contraction
\bea s\Big(f_0[f_1|\cdots|f_q]m\Big)=1[f_0|f_1|\cdots|f_q]m,~~~\{b_1,s\}=\textrm{id}-\eta\epsilon.\label{contraction}\eea

The bar complex has a second grading, which assigns to an element $f_0[f_1|\cdots|f_q]m$
\bea \deg f_0+\sum_{i=1}^q(\deg f_i-1)+\deg m.\label{grading}\eea
There is another internal differential $b_0$ coming from the differential of \cdga~and \dgm~structures (denoted $Q$ and $\hat Q$ respectively),
\bea b_0\Big(f_0[f_1|\cdots|f_q]m\Big)&=&(Qf_0)[f_1|\cdots|f_q]m\nn\\
&&-\sum_{i=1}^{q}f_0[f_1|\cdots |f_{i-1}|Qf_i|f_{i+1}|\cdots|f_q]m(-1)^{f_0+(f_1+1)+\cdots +(f_{i-1}+1)}\nn\\
&&+f_0[f_1|\cdots|f_q](\hat Qm)(-1)^{f_0+(f_1+1)+\cdots +(f_q+1)}.\label{diff_internal}\eea
Since $\{b_0,s\}=0$ and from Eq.\ref{contraction} we have
\bea H^{\sbullet}_b(B(\A,\A,\Mod))=H^{\sbullet}_{\hat Q}(\Mod)\nn\eea
and $B(\A,\A,\Mod)$ is a resolution of $\Mod$. One can also define a more general bar complex by the tensoring over $\A$
\bea B^{\sbullet}(\Nod,\A,\Mod)=\Nod\otimes_{\A}B^{\sbullet}(\A,\A,\Mod),\label{def_bar_gen}\eea
where $\Nod$ is a (right) \dgm~ over $\A$.

It is common to normalize the bar complex so as to get rid of negative grading elements (c.f. Eq.\ref{grading}), which is indispensable for the use of spectral sequence later. In the case $\A$ is \emph{connected}, namely, the degree 0 part of $\A$ is $k$ itself, the normalized bar complex is defined
\bea \bar B^{-q}(\A,\A,\Mod)=\A\otimes_k\underbrace{\A^+\otimes_k\A^+\cdots\A^+}_q\otimes_k \Mod,\label{def_bar}\eea
where $\A^+$ are those $x\in\A,\,\deg x>0$ and the resolution $\bar B^{\sbullet}$ is proper projective (\cite{McCleary} ch.7). The normalized bar complex $\bar B^{\sbullet}(\Nod,\A,\Mod)=\Nod\otimes_{\A}\bar B^{\sbullet}(\A,\A,\Mod)$ is defined similarly. The next is a classic result of MacLane,
\begin{Lem}\label{Lem_Maclane}(\cite{Maclane} ch.8.6), Assume $\A^0=k$, then
\bea H^{\sbullet}_{b_1}(B(\Nod,\A,\Mod))=H^{\sbullet}_{b_1}(\bar B(\Nod,\A,\Mod)).\label{Maclane}\eea
\end{Lem}
This lemma actually applies to any simplicial module. In the case $\Nod=\A$, this can be seen directly by using a slightly modified contraction
\bea s\Big(f_0[f_1|\cdots|f_q]m\Big)=\bigg\{\begin{array}{c}
                                              1[f_0|f_1|\cdots|f_q]m\\
                                              0
                                            \end{array}
                                            \begin{array}{c}
                                              \textrm{if }\deg f_0>0  \\
                                              ~~\textrm{otherwise}
                                            \end{array},\nn\eea
we can check that this is still a contraction provided $\A^0=k$, leading to the acyclicity.

In the case $\A$ is not connected which is relevant for $\A=\Omega^{\sbullet}(M)$, a slight modification to the normalization is needed. Following the notation of ref.\cite{cyc_bar} we define the operator $S_i(f)$, where $f\in\A^0$, that acts on $B^{-q}(\Nod,\A,\Mod)$ as
\bea S_i(f)\Big(n[f_1|\cdots|f_q]m\Big)=f_0[f_1|\cdots|f_{i-1}|f|f_i|\cdots f_q]m,~~~1\leq i\leq q+1.\nn\eea
Let the normalized bar complex $\bar B^{-q}$ be a quotient of $B^{-q}$
\bea \bar B^{-q}=B^{-q} /(\textrm{img}\,S_i(f),\,\textrm{img}\,\{S_i(f),b\}).\label{quotient}\eea
The quotient by $\textrm{img}\,S_i(f)$ restricts the $\A$'s 'between the bars' to have positive degrees, while further modding out by $\textrm{img}\,\{S_i(f),b\}$ gives us relations like
\bea &&-n[f_1|Qf|f_2]m+n[f_1f|f_2]m-n[f_1|ff_2]m\sim0.\label{relations}\eea
that are eventually related to the quotient by eom in Eq.\ref{ideal}.

We first note that when $\Nod=\A$, then the contraction Eq.\ref{contraction}
pass over to the quotient, thus we still have
\bea H^{\sbullet}_b(\bar B(\A,\A,\Mod))=H^{\sbullet}_{\hat Q}(\Mod).\label{still_true}\eea
For general $\Nod$, the cohomology is altered after the normalization. However, if one assumes that $H^0_Q(\A)=k$, one can prove the
\begin{Lem}(see ref.\cite{cyc_bar} Lem.5.3 and the remark thereafter)
If $H^0(\A)=k$, then
\bea H_b^{\sbullet}(\bar B(\Nod,\A,\Mod))=H_b^{\sbullet}(B(\Nod,\A,\Mod)).\nn\eea
\end{Lem}
The proof uses explicit representatives for $\bar B^{-q}(\Nod,\A,\Mod)$ instead of working with the quotient Eq.\ref{quotient}. Indeed
$\bar B^{\sbullet}$ can be described as
\bea \bar B^{-q}(\Nod,\A,\Mod)\sim~\Nod\otimes_k\underbrace{\bar\A\otimes_k\bar\A\cdots\bar\A}_q\otimes_k \Mod.\label{explicit}\eea
where $\bar \A$ is a subalgebra of $\A$ and $(Q\A^0)^{\perp}$ is a complement of $Q\A^0$ in $\A^1$
\bea \bar\A^q=\Bigg\{\begin{array}{c}
            0\\
            (Q\A^0)^{\perp}\\
            \A^q
          \end{array}\begin{array}{c}
                       q=0 \\
                       q=1 \\
                       q>1
                     \end{array}.\nn\eea
Proof: There is a map $\bar B^{\sbullet}\to B^{\sbullet}$ is induced by $\bar\A\to\A$.
Next one filtrates $\bar B^{-q}(\Nod,\A,\Mod)$ and $B^{-q}(\Nod,\A,\Mod)$ according to $q$,
\bea F^{-q}B=\Large\textrm{$\oplus$}_{i=0}^{i=q}B^{-i},~~~~\bar F^{-q}\bar B=\Large\textrm{$\oplus$}_{i=0}^{i=q}\bar B^{-i}\nn\eea
The $E_1$ term for both is
\bea E_1^{-p,\sbullet}=B^{-p}(H(\Nod),H(\A),H(\Mod)),~~~\bar E_1^{-p,\sbullet}=B^{-p}(H(\Nod),H(\bar\A),H(\Mod))\nn\eea
and their $E_2$-terms are the $b_1$ cohomology of $E_1$. Since $H(\bar \A)$ and $H(\A)$ differ only by $k$ at degree 0, thus $\bar E^{-p,q}$ is the normalized version of $E^{-p,q}$ for every $q$, and that they have the same $b_1$ cohomology follows from Lem.\ref{Lem_Maclane}.
The isomorphism at $E_2$-term shows that there is an isomorphism of the entire bar complex$\square$

For our interest, we consider a representation, a $Q$-equivariant bundle $(\E,\hat Q)\stackrel{\pi}{\to}(\Ta,Q)$. Take now
$\A=C^{\infty}(\Ta)$ and $\Mod=\Gamma(\E)$, then through the map $\pi$, the module $\Mod$ is a \dgm~over $\A$, or we may of course take $\Mod$ to be $\A$ itself.

There are two variants to the above definition of bar complex that we are also interested in. First the fix-ended bar complex: $\bar B^{-q}(k,\A,k)$, where $k$ receives an $\A$-module structure through an augmentation $\A\stackrel{\epsilon}{\to}k$. The typical elements are denoted $[f_1|\cdots|f_q]$, the differential $b_0$ remains the same. As for $b_1$, since we can take $f_i$ to be of positive degree, $\epsilon f_i=0$, eliminating all but the second line of Eq.\ref{diff_bar} and leaving
\bea
b_0\Big([f_1|\cdots|f_q]\Big)&=&-\sum_{i=1}^{q}[f_1|\cdots |f_{i-1}|Qf_i|f_{i+1}|\cdots|f_q](-1)^{(f_1+1)+\cdots +(f_{i-1}+1)}\nn\\
b_1\Big(\left[f_1|f_2|\cdots|f_q\right]\Big)&=&\sum_{i=1}^{q-1}[f_1|\cdots |f_{i-1}|f_if_{i+1}|f_{i+2}|\cdots|f_q](-1)^{(f_1+1)+\cdots +(f_i+1)}.\label{diff_bar_based}\eea

Another variant to the bar complex is known as the cyclic-bar complex \cite{cyc_bar,GJPbar}, which is denoted as $C^{-q}(\A)$
\bea C^{-q}(\A)=\A\otimes_k\underbrace{\A\otimes_k\A\cdots\A}_q,\label{def_bar_cyc}\eea
and a typical element is denoted
\bea f_0[f_1|\cdots|f_q]\in C^{-q}(\A),\label{typical_cyc_bar}\eea
on which the differential $b_1$ acts as
\bea   b_1\Big(f_0\left[f_1|f_2|\cdots|f_q\right]\Big)&=&(-1)^{f_0}f_0f_1[f_2|f_3|\cdots|f_q]\nn\\
&&+\sum_{i=1}^{q-1}f_0[f_1|\cdots |f_{i-1}|f_if_{i+1}|f_{i+2}|\cdots|f_q](-1)^{f_0+(f_1+1)+\cdots +(f_i+1)}.\nn\\
&&-f_qf_0[f_1|\cdots|f_{q-1}].\label{diff_bar_cyclic}\eea
The normalized version $\bar C^{\sbullet}$ is defined similarly. This bar complex is used in ref.\cite{2010arXiv1006.1240Q} to construct $L_{\infty}$-algebra weight systems for knots.

\subsection{Iterated Integral}
Let $\So=T[1][0,1]$ or $\So=T[1]S^1$, one can map the bar complex to functions on the space of \cdga-maps $\So\stackrel{\varphi}{\to} \Ta$. The construction is through the iterated integral, similar to the Wilson-line operators in physics. We first see how this is done for $\Ta=T[1]M$. In this case, the functions $C^{\infty}(\Ta)$ are differential forms on $M$, and space of \cdga-maps $\So\stackrel{\varphi}{\to}T[1]M$ is the path (loop) space $PM~(LM)$. Forms $\Omega^{\sbullet}(M)$ can be pulled back into forms $\Omega^{\sbullet}(PM)\otimes \Omega^{\sbullet}([0,1])$, by simply writing
\bea dx^{\mu}\to \dot \varphi^{\mu}(t)dt+ \delta \varphi^{\mu}(t),\label{pull_back_easy}\eea
where $\varphi$ is a map $[0,1]\to M$ and $\delta\varphi$ is the variation of the map $\varphi$, taken as a 1-form on $PM$.

Denote by $\Delta^q$ standard $q$-simplex parameterized as $\Delta^q=\{(t_1,\cdots t_q)|1\geq t_1\geq \cdots\geq t_q\geq0\}$. Then
for a collection of forms $f_1,\cdots f_q$ on $M$, we have a series of maps
\bea f_1\otimes f_2\cdots \otimes f_q\longrightarrow \Omega^{\sbullet}(PM)\otimes\Omega^{\sbullet}(\Delta^q)\stackrel{\int_{\Delta^q}}{\longrightarrow }\Omega^{\sbullet}(PM),\nn\eea
where in the first step, one applies Eq.\ref{pull_back_easy}, and the second map is the integration over $\Delta^q$.
%
In the special case when all $f_i$ are 1-forms, the above map is simply written as the iterated integral
\bea f_1\otimes f_2\cdots \otimes f_q\to \int^1_0 \varphi^*f_1(t_1) \int^{t_1}_0 \varphi^*f_2(t_2)\cdots \int^{t_{q-1}}_0 \varphi^*f_q(t_q)\in \Omega^0(PM),\nn\eea
and the iterated integral gives a function on $PM$. Furthermore, if the 1-forms $f$ are a connection, then this expression is just (part of) the Wilson-line operator. But in general one gets higher forms in $PM$.

It is worthwhile repeating the above calculation from the point of view of sec.\ref{sec_SF}. Consider the mapping problem
$\So\times\Pa \to T[1]M$, with an arbitrary choice of a test manifold $\Pa$. The superfields are
\bea \BS{x}^{\mu}(t,\theta,y),~~~\BS{v}^{\mu}(t,\theta,y),\nn\eea
where $y$ is the generator of $C^{\infty}(\Pa)$. Write out Eq.\ref{condition_dga},
\bea &&0=\BS{v}^{\mu}(t,\theta,y)-(\theta\partial_t+d^y)\BS{x}^{\mu}(t,\theta,y)\nn\\
&&0=-(\theta\partial_t+d^y)\BS{v}^{\mu}(t,\theta,y),\nn\eea
and expand in power of $\theta$
\bea 0=v^{\mu}_t(t,y_0)-\partial_tx^{\mu}(t,y_0),~~0=v^{\mu}(t,y)-d^yx^{\mu}(t,y_0),~~0=-\partial_tv^{\mu}(t,y)+d^yv_t^{\mu}(t,y_0).\nn\eea
The second equation again says that $v^{\mu}(t,y)$ is regarded as a 1-form in $\Omega^{\sbullet}(PM)$, compare with Eq.\ref{pull_back_easy}.

Letting now $f_i$ be functions $C^{\infty}(T[1]M)$, $z$ be the collective coordinate $(t,\theta)$ and suppress completely the coordinate $y$ of the test manifold, then the iterated integral is written in super language as
\bea &&f_0[f_1| \cdots |f_q]f_{q+1}~\to~f_0(1)\smalint_0^1 dz_1 \BS{f}_1(z_1)\smalint^{t_1}_0dz_1\BS{f}_2(z_2)\smalint^{t_2}_0\cdots\smalint^{t_{q-1}}_0dz_q\BS{f}_q(z_q)\cdot f_{q+1}(0)\nn\\
&&\hspace{1cm}=f_0(1)\int_0^1 dt_1\Big( [\partial_{\theta_1}\BS{f}_1]\Big|_{\theta_1=0}\Big)\int_0^{t_1} dt_2\Big( [\partial_{\theta_1}\BS{f}_2]\Big|_{\theta_1=0}\Big)\cdots\int^{t_{q-1}}_0dt_q\Big([\partial_{\theta^q}\BS{f}_q]\Big|_{\theta^q=0}\Big)
f_{q+1}(0)\label{int_iter},\eea
where we write $\BS{f}_0(z)\big|_{t=1,\theta=0}$ simply as $f_0(1)$, and $\BS{f}_{q+1}(z)\big|_{t=0,\theta=0}$ as $f_{q+1}(0)$.

This iterated integral is generalized verbatim to mapping problems $T[1][0,1]\otimes \Pa\stackrel{\varphi}{\to}\Ta$ with arbitrary $NQ$-manifold $\Ta$.
\begin{rmk}\label{rmk_neg_deg}
Also recall that in a superfield expansion, any component of negative degree are the auxiliary fields are set to zero at the end of the calculation, thus if any $f_i$ in the above expression is of degree 0 then the whole integral vanishes. Here we see why the entries 'between the bars' in $B^{-q}$ have to belong to $\A^+$ (c.f. the quotient by $\textrm{img}\,S_i(f)$ in Eq.\ref{quotient}).
\end{rmk}

We have seen that for $\M=T[1]M$ the image of the iterated integral is the differential forms $\Omega^{\sbullet}(PM)$, equivalently one may say that
$T[1]PM$ is the $NQ$-manifold representing the presheaf $\HomM_{NQ}(T[1][0,1],T[1]M)$. For a general $\Ta$, we can at best say the image of Eq.\ref{int_iter} is a function on the space of super-curves $(C^{\infty}(\bar{\cal S}),\gd_B)$.

%
%
%
%
%

Denoting $\A=C^{\infty}(\M)$, we have an important lemma
\begin{Lem}
Under the iterated integral, the total differential $b_0+b_1$ of the bar complex is sent to the differential $\delta_B$,
\bea (\bar B^{\sbullet}(\A,\A,\A),b_0+b_1)\stackrel{\textrm{iterated integral}}{\longrightarrow}(C^{\infty}(\bar{\cal S}),\gd_B).\label{bar_inter_it}\eea
\end{Lem}
The proof is largely taken from ref.\cite{GJPbar}, yet here we redo the calculation in terms of superfields and also for the reason that the same manipulation will be used over and over again in showing the properties of Wilson lines later.

\noindent proof (sketch):
Consider $q$ functions on $NQ$-manifold $(\M,Q)$, and pick a \cdga~map $T[1][0,1]\times\Pa\to\M$, and we will do all the calculation first in ${\cal S}$ and then pass onto $\bar{\cal S}$. Apply $\gd_B$ to the rhs of Eq.\ref{int_iter}
\bea &&\delta_B\Big(f_0(1)\smalint^{1}_0dz_1\BS{f}_1(z_1)\smalint^{t_1}_0\cdots\smalint^{t_{q-1}}_0dz_q\BS{f}_q(z_q)\cdot f_{q+1}(0)\Big)\nn\\
&&\hspace{2cm}=(Qf_0(1))\smalint^{1}_0dz_1\BS{f}_1(z_1)\smalint^{t_1}_0\cdots\smalint^{t_{q-1}}_0dz_q\BS{f}_q(z_q)\cdot f_{q+1}(0)\nn\\
&&\hspace{2.5cm}+\sum_{i=1}(-1)^{f_0+(f_1+1)+\cdots+(f_{i-1}+1)}\smalint_0^1 dz_1 \BS{f}_1(z_1)\cdots \smalint^{t_{i-2}}_0dz_{i-1}\BS{f}_{i-1}(z_{i-1})\nn\\
&&\hspace{4.5cm}\int^{t_{i-1}}_0dt_i\partial_{\theta_i}\big[D\BS{f}_i(z_i)-\BS{Q}\cdot \BS{f}_i(z_i)\big]\Big|_{\theta_i=0}\smalint_0^{t_i}
               \cdots\smalint^{t_{n-1}}_0dz_q\BS{f}_q(z_q)\cdot f_{q+1}(0)\nn\\
&&\hspace{2.5cm}+(-1)^{f_0+(f_1+1)+\cdots+(f_q+1)}f_0(1)\smalint^{1}_0dz_1\BS{f}_1(z_1)\smalint^{t_1}_0\cdots\smalint^{t_{q-1}}_0dz_q\BS{f}_q(z_q)\cdot Qf_{q+1}(0).\label{computation}\eea
The first, last and the second term in the square brace in the middle term produces the the three lines of Eq.\ref{diff_internal}.
The first term in the square brace gives
\bea \cdots\int^{t_{i-1}}_0dt_i\partial_{\theta_i}\big[D\BS{f}_i(z_i)\big]\Big|_{\theta_i=0}\int_0^{t_i}
               \cdots&=&\cdots\int^{t_{i-1}}_0dt_i\partial_{t_i}f_i(t_i)\int_0^{t_i}\cdots\nn\eea
An integration by part will either collapse two neighboring $f_i$ and $f_{i+1}$ together, giving the second line of Eq.\ref{diff_bar}; or collapse $f_1$ with $f_0$ or $f_q$ with $f_{q+1}$, giving the first and third line of Eq.\ref{diff_bar}. Finally one may pass from ${\cal S}$ to $\bar{\cal S}$ by taking the quotient Eq.\ref{ideal}, leading to Eq.\ref{bar_inter_it}$\square$

Note that taking this quotient parallels the quotient Eq.\ref{quotient} on the bar complex side, we leave it to the reader to show this himself. In short, the superfield technique plus keeping the auxilliary fields permits us to perform computations 'off-shell'.

%

\subsection{Variants to the Iterated Integral}
In sec.\ref{BCaII}, two variants to the bar complex were introduced, the iterated integral for these two types are slightly modified.

Let $\bar{\cal S}_{x_1,x_0}$ denote the space of super-curves with ends fixed at $x_1,x_0\in|\M|$.
Take one such curve $T[1][0,1]\times \Pa\to \Ta$, it is parameterized as $x\to \BS{x}(t,\theta,y)$ such that $\BS{x}^m(1,\theta,y)\big|_{\theta=0}=\BS{x}^m(0,\theta,y)\big|_{\theta=0}=0$ if $\deg x^m>0$ and $\BS{x}^m(1,\theta,y)=x^m_1$, $\BS{x}^m(0,\theta,y)=x^m_0$ if $\deg x^m=0$. Since $f_i\in \A^+,~i=1\cdots q$, we have clearly $f_1(1)=f_q(0)=0$ due to the boundary condition, then it should be clear that the iterated integral gives a homomorphism
\bea (\bar B^{\sbullet}(k,\A,k),b_0+b_1)\stackrel{\textrm{iterated integral}}{\longrightarrow}(C^{\infty}(\bar{\cal S}_{x_1,x_0}),\gd_B),\label{bar_inter_it_fixed}\eea
where the $\A$-module structure of the first $k$ is induced from the evaluation map at $x_1$ and the second induced by the evaluation map at $x_0$.

Next we have a lemma that basically says the space of paths in $\M$ has the same cohomology as $\M$ itself.
\begin{lem}\label{lem_shrink}
$H^{\sbullet}_{\delta_B}(\bar{\cal S})=H_Q^{\sbullet}(\M)$
\end{lem}
proof: Roughly the argument is that every curve can be shrunk to its initial point giving the homotopy equivalence of $\bar{\cal S}$ and $\M$. But since now the notion of homotopy is a little different from that of the ordinary path space, we run the argument again here.

If $\BS{x}^m(t,\theta,y)$ parameterizes the super-curve, then the homotopy $T[1][0,1]\times \Pa\times T[1][0,1]\to \M$ parameterized as
\bea \BS{X}^m(t,\theta,s,\vartheta,y)=\BS{x}^m(ts,s\theta+t\vartheta,y)\nn\eea
shrinks the super curve to the constant curve which is the initial point.

Let $\omega$ be a $\gd_B$-closed class in $\bar{\cal S}$, thus for any test manifold, $d^y\omega(y)=0$ since under the mapping $\Pa\to \bar{\cal S}$, $d^y$ goes to $\gd_B$. We have
\bea \omega(y)\big|_{s=1}-\omega(y)\big|_{s=0}=\int_0^1 ds~\frac{\partial}{\partial s}\omega(y)=\big(\{\gd_B,\iota_{\bar x}\}+d^y\bar x^m\frac{\partial}{\partial x^m}\big)\int_0^1 ds \omega(y),\nn\eea
using $(\gd_B\BS{x})(y)=d^y\BS{x}(y)$, the above simplifies to
\bea \omega(y)\big|_{s=1}-\omega(y)\big|_{s=0}=d^y\int_0^1 ds~\iota_{\bar x}\omega(y).\nn\eea
If we note that at $s=0,\vartheta=0$, $\gd_B$ acts as $Q$, then the result follows$\square$

The next proposition is a generalization of a similar statement in ref.\cite{GJPbar}, which shows that the image of $\bar B^{\sbullet}(k,\Omega^{\sbullet}(M),k)$ under the iterated integral computes the cohomology of the space of path with fixed ends. Unfortunately, the following generalization to $NQ$-manifolds requires some excessively strong conditions on $\M$.

\begin{Prop}\label{iso_bar_fixed}Assume $\M$ is path connected, simply connected and $H_Q^1(\M)=0$, then the iterated integral
\bea \bar B^{\sbullet}(k,\A,k)\to C^{\infty}(\bar{\cal S}_{x_1,x_0})\nn\eea
induces isomorphism in cohomology $H_{b}(\bar B(k,\A,k))\to H_{\gd_B}(\bar {\cal S}_{x_1,x_0})$, where $\bar{\cal S}_{x_1,x_0}$ is the space of super-curves with ends fixed.
\end{Prop}
\begin{rmk}
The statement '$\M$ is path connected' is in the sense that a \dga~map $T[1][0,1]\to\M$ exists for any initial, final points $x_0,x_1\in|\M|$. If, for example, $\M$ is as in Ex.\ref{Lie_algebroid}, then the surjectivity of the anchor $\rho$ implies the statement. And in general, path connectedness implies $H^0_Q(\M)=H^0_Q(\A)=k$; the converse is not true. It is also unclear to us if connectedness plus simply-connectedness implies automatically $H_{Q}^1(\M)=H^1_{Q}(\A)=0$.

The first two assumptions are required in order to canonically identify the groups $H^{\sbullet}_{\gd_B}(\bar {\cal S}_{x_1,x_0})$ for different $x_{1,0}$. While $H_Q^1(\A)=0$ is needed to show $H^0_b(\bar B(k,\A,k))=k$.\end{rmk}
proof:
One filtrates both $\bar B(\A,\A,\A)$ and $C^{\infty}(\bar{\cal S})$ and computes their cohomology using a spectral sequence argument, while the iterated integral is a morphism between the two spectral sequences.

The filtration on $\bar B(\A,\A,\A)$ is that we say $f_0[f_1|\cdots|f_q]f_{q+1}$ belongs to filtration $s$ if $\deg f_0+\deg f_{q+1}\geq s$; this filtration is clearly preserved by both $b_0,\,b_1$. The $E_1$-term of this filtration is
\bea E^{r,s}_1=H_b^{r+s}(F^r/F^{r+1}\bar B),\nn\eea
by using the explicit description of the reduced bar complex Eq.\ref{explicit} plus the filtration condition, we have
\bea E^{r,s}_1=\bigoplus_{p+q=r}\A^p\otimes_kH_b^{r+s}(\bar B(k,\A,k))\otimes_k\A^q,\nn\eea
and hence $E_2$ term is
\bea E_2^{r,s}=H^r_Q(\A\otimes_k\A)\otimes H^s_b(\bar B(k,\A,k))\simeq E_2^{r,0}\otimes E_2^{0,s},\nn\eea
where the last isomorphism uses $H^0_Q(\A\otimes_k\A)=H^0_b(\bar B(k,\A,k))=k$. We also know that $E_{\infty}=H_b^{\sbullet}(\bar B(\A,\A,\A))=H^{\sbullet}_Q(\A)$ from Eq.\ref{still_true}, and $H^{\sbullet}_Q(\A)$ is just what we call $H^{\sbullet}_Q(\M)$.

We also consider the Serre fibration
\bea
\begin{tikzpicture}
  \matrix (m) [matrix of math nodes, row sep=2.5em, column sep=.5em]
    {  \bar{\cal S} & \bar{\cal S}_{x_1,x_0}\\
       \textrm{\scriptsize{$\widehat{\M\times \M}$}} &  \\ };
 \path[->]
 (m-1-2) edge  (m-1-1)
 (m-1-1) edge  node[right] {\textrm{\scriptsize{$t\times h$}}} (m-2-1) ;
\end{tikzpicture}\nn\eea
where $t\times h$ projects a super-curve to the final and initial point. We say a function is in filtration $F^sC^{\infty}(\bar{\cal S})$ if it is in the ideal generated by $(t\times h)^* \omega$ with $\omega$ being all the functions on the base with degree $\geq s$. Using this filtration, we have the $E_2$ term \bea \bar E_2^{r,s}=H^r_Q(\M\times \M)\otimes H_{\gd_{B}}^s(\bar{\cal S}_{x_1,x_0})\simeq \bar E_2^{r,0}\otimes \bar E_2^{0,s},\nn\eea
where the last step requires $H^0(\M)=k$ and $H^0_{\gd_B}(\bar {\cal S}_{x_1,x_0})=k$.

We know from Lemma.\ref{lem_shrink} that $\bar E_{\infty}=H^{\sbullet}_{Q}(\M)$, thus the iterated integral induces an isomorphism $E_{\infty}\simeq \bar E_{\infty}$. On the other hand it is clear the iterated integral also induces an ismorphism $E^{r,0}_2\simeq \bar E^{r,0}_2$, thus we appeal to Zeeman's comparison theorem (\cite{McCleary} ch.3) and conclude $E^{0,s}_2\simeq \bar E^{0,s}_2$ under the iterated integral. This exactly says the image of $\bar B(k,\A,k)$ under the iterated integral computes the cohomology of space of super-curves with fixed ends $\square$

The second variation of the bar complex is the cyclic bar complex. The corresponding iterated integral is a bit lengthy to write, the reader may see ref.\cite{2010arXiv1006.1240Q} Eq.24 or ref.\cite{GJPbar} Thm 2.1, the formula roughly corresponds to picking an arbitrary base point on the loop as time 1, writing the iterated integral as before and finally integrating over the position of the base point along the loop.

Using a filtration on the degree of $f_0$ for an element $f_0[f_1|\cdots|f_q]$ in complex $\bar C^{-q}(\A)$, and a filtration according to the base degree in the fibration (where $\bar{\cal S}_{x,x}$ is the space of based super loops, $\bar{\cal S}_{\textrm{O}}$ unbased loops)
\bea
\begin{tikzpicture}
  \matrix (m) [matrix of math nodes, row sep=2em, column sep=2em]
    {  \bar{\cal S}_{\textrm{O}} &{\bar{\cal S}}_{x,x} \\
       \widehat{\M} &  \\ };
 \path[->]
 (m-1-2) edge  (m-1-1)
 (m-1-1) edge  node[right] {} (m-2-1) ;
\end{tikzpicture}\nn\eea
one can show similarly
\begin{Prop}\label{iso_bar_loop}Under the assumption of prop.\ref{iso_bar_fixed},
the iterated integral induces a map $\bar C^{\sbullet}(\A)\to C^{\infty}(\bar{\cal S}_{\textrm{O}})$, which is an isomorphism in cohomology .
\end{Prop}


\section{Wilson Lines}\label{sec_WL}
In this bulky section, we will apply the iterated integral to a particular member of the bar complex and obtain the formula of the Wilson line. We shall first work locally in a trivialization patch and then connect pieces of Wilson lines in different patches with a transition function and show that the Wilson line do not in general transform covariantly; we also work out how the Wilson-lines change under homotopies.

Pick a representation of the $Q$ structure on $\M$, a $Q$-equivariant bundle $(\E,\hat Q)$ covering $(\M,Q)$, and $\Mod=\Gamma(\E)$ becomes a \dgm~over the (not necessarily commutative) \dga~$\Gamma(\textrm{End}\,\E)$. We can form the bar complex $\bar B^{\sbullet}(\Mod,\Gamma(\textrm{End}\,\E),\Mod^*)$ as in sec.\ref{BCaII}, but we shall \emph{work locally first} and treat members of $\bar B^{\sbullet}(\Mod,\Gamma(\textrm{End}\,\E),\Mod^*)$ as members of $\bar B^{\sbullet}(\Mod,C^{\infty}(\M),\Mod^*)$ under a trivialization.

Ideally, an operator for the parallel transport should be a section of the pull-back bundle $t^*\E\otimes h^*\E^*$ over $\bar{\cal S}$ which we denote as $\E\textrm{nd}$, where $h,\,t$ are the projections
\bea
&&\E\textrm{nd}\to\E\;\times\;\E^*\nn\\
&&\,\downarrow\hspace{.9cm}\downarrow~~~~\downarrow\nn\\
&&\bar{\cal S}\stackrel{t\times h}{\longrightarrow}\M\times\M\label{ideal_situation}\eea
of a super curve to its first and last point. The differential $\gd_B$ on $\bar{\cal S}$ has a natural lift to $\E\textrm{nd}$ just like $Q$ on $\M$ is lifted to $\hat Q$ on $\E$. As such, a Wilson line should be a presheaf sending $\Pa$ to the set of \dgm~maps $\Hom_{\dgm}(\Gamma(\E\textrm{nd}),C^{\infty}(\Pa))$ that cover the \cdga~maps $\Hom_{\cdga}(\Gamma(C^{\infty}(\bar{\cal S}),C^{\infty}(\Pa))$, were it not for trivialization problems.

\subsection{Local Definition of Wilson-Line}
We choose the fibre coordinate of $\E$ as $\zeta_{\ga}$, in other words, $\zeta_{\ga}$ is the local generators of the $C^{\infty}(\M)$-module $\Gamma(\E)$. We also let $\zeta^{\ga}$ be the dual fibre coordinate of $\E^*$ (pay attention to the position of the index). Let the representation matrix $T^{\ga}_{~\gb}$ be defined in sec.\ref{RoQS}. Locally we can form an element of the bar complex
\bea W_q=\zeta_{\ga}(W_p)^{\ga}_{~\gb}\zeta^{\gb}\in B^{-q}(\Mod,\A,\Mod^*),~~~(W_q)^{\ga}_{~\gb}=[T^{\ga}_{~\gc_1}| T^{\gc_1}_{~\gc_2}|\cdots|T^{\gc_{q-1}}_{~~\gb}],\,q>0,~~~~(W_0)^{\ga}_{~\gb}=\delta^{\ga}_{~\gb}.\label{vissnamn}\eea
Out of $W_q$ we can make a formal series
\bea W=\sum_{q=0}^{\infty}~W_q.\label{vissnamn1}\eea

We also pick a super-curve $\BS{C}:~T[1][0,1]\times\Pa\to \M$. Though this curve may cross several coordinate patches, we may subdivide it into pieces and thus without loss of generality, we may assume $T[1][0,1]\times\Pa\to \M$ is in one patch only.

Then in one coordinate patch, the Wilson line is defined as in
\begin{defi}\label{def_W}For a given super-curve $C$,
the Wilson-line is the image of $W$ in Eq.\ref{vissnamn1} under the iterated integral
\bea
W\stackrel{\textrm{iterated integral}}{\longrightarrow}~U_C=\zeta_{\ga}U^{\ga}_{C\gb}(1,0)\zeta^{\gb},\nn\eea
concretely, $U^{\alpha}_{C\beta}(1,0)$ is given by the familiar path-ordered integral~\footnote{We remind the reader that the boldface symbols are promoted to being superfields, e.g. $\BS{T}=T(\BS{x}(t,\theta))$}
\bea U^{\alpha}_{C\beta}(1,0)=\left(\BB{P}\exp\Big(-\int_0^1dtd\theta~\BS{T}\Big)\right)^{\alpha}_{~\beta}\label{wilson_line},\eea
where $\BB{P}$ is the path ordering operator which has the same effect as the iterated integral.\end{defi}
The subscript $C$ will usually be omitted, and $\zeta_{\ga}\,(\zeta^{\gb})$ on the rhs is understood as the coordinate of $\E\,(\E^*)$ over the final (initial) point of $C$.

\begin{figure}[h]
\begin{center}
\begin{tikzpicture}[scale=.8]
\draw [blue] (0,0) circle (1);
\draw [blue] (1.2,.4) circle (1);
\draw (-.5,0) arc (270:315:3);
\node (a) at (.6480,0.228) [fill,circle,inner sep=1pt] {};
\node (b) at (-.5,0) [fill,circle,inner sep=0pt] {};
\node (c) at (1.633,0.854) [fill,circle,inner sep=0pt] {};
\draw (a) node [above] {\small$t$}
 (b) node [above] {\scriptsize{0}}
 (c) node [above] {\scriptsize{1}};
\end{tikzpicture}\caption{Two pieces of Wilson-lines connected by a transition function}\label{splicing}
\end{center}
\end{figure}
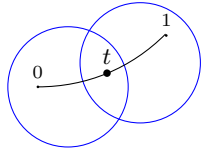

In the case the Wilson-line/loop passes through several patches each with its own trivialization, as in fig.\ref{splicing}.
Then one splices together Wilson-lines by using transition functions; it suffices to look at the case of two patches, let $\Omega^{\gc}_{~\tilde\gb}$ be the transition function, then the spliced Wilson-line is given as
\bea \zeta_{\ga}U^{\ga}_{~\tilde\gb}(1,0)\tilde\zeta^{\tilde \gb}=\zeta_{\ga}U^{\ga}_{~\gc}(1,t)\cdot\Omega^{\gc}_{~\tilde\gd}(t)\cdot U^{\tilde \gd}_{~\tilde\gb}(t,0)\tilde\zeta^{\tilde \gb},\nn\eea
where $\Omega^{\gc}_{~\tilde\gd}(t)=\Omega^{\gc}_{~\tilde\gd}(\BS{x}(t,\theta,y))\big|_{\theta=0}$ is the transition function 'at the connecting point'. The ability to cut the Wilson-line into pieces and reconnect them requires $U^{\ga}_{~\gb}(1,0)$ to be invariant under reparameterization, which will be proved in cor.\ref{cor_repara}.

\begin{rmk}
By the phrase 'Wilson line' we mean either those written locally in one patch or those spliced together, as the case may be. Indeed, one will see later that the problem of splicing is rather trivial.
\end{rmk}

\begin{rmk}
A simple counting shows that the level $\ga$ must be no greater than $\gb$, because $\deg U^{\ga}_{~\gb}+\ga-\gb=0$ and $\deg U^{\ga}_{~\gb}$ must be no less than 0.
\end{rmk}

At least locally, we can regard the Wilson-line as a function of $C^{\infty}(\bar{\cal S})$ with values in $\Hom(\E|_0,\E|_1)$. The BRST differential $\delta_B$, which is the image of the differential $b$ in Eqs.\ref{diff_bar},\ref{diff_internal}, acts on $C^{\infty}(\bar{\cal S})$, $\zeta_{\ga}$ and $\zeta^{\gb}$ as
\bea &&\delta_B\BS{x}^A(t,\theta,y)=-D\BS{x}^A(t,\theta,y)+\BS{Q}^A(t,\theta,y),\nn\\
&&\delta_B\zeta_{\ga}=\hat Q\zeta_{\ga}=(-1)^{\gb}\zeta_{\gb}T^{\gb}_{~\ga},\nn\\
&&\delta_B\zeta^{\ga}=\hat Q\zeta^{\ga}=-T^{\ga}_{~\gb}\zeta^{\gb}.\label{BRST}\eea
\begin{Prop}\label{FIP}
The BRST differential acts on the Wilson line as
\bea \delta_B U^{\alpha}_{~\beta}(1,0)=-T^{\alpha}_{~\delta}(1)U^{\delta}_{~\beta}(1,0)+(-1)^{\alpha+\gamma}U^{\alpha}_{~\gamma}(1,0)T^{\gamma}_{~\beta}(0)
\label{BRST_WL},\eea
where $T(0)=T(\BS{x}(0,\theta,y))\big|_{\theta=0}$, $T(1)=T(\BS{x}(1,\theta,y))\big|_{\theta=0}$ act on the fibre of $\E$ over the initial and final point of the path.

By including the action of $\gd_{B}$ on $\zeta_{\ga},\,\zeta^{\ga}$, one can equally write Eq.\ref{BRST_WL} concisely as
\bea \delta_BU_C=0.\label{concise}\eea
\end{Prop}

\begin{defi}\label{def_W_loop}For a super-loop we define the Wilson loop to be
\bea U_{\textrm{O}}=\sum_{\alpha}~U^{\alpha}_{~\alpha}(1,0)(-1)^{\alpha}.\label{Wilson_loop}\eea
\end{defi}
An important corollary is
\begin{Cor}
The Wilson loop is $\delta_B$ closed.
\end{Cor}
Proof:
\bea&&\hspace{1cm}\delta_BU_{\textrm{O}}=-T^{\alpha}_{~\delta}U^{\delta}_{~\alpha}(-1)^{\alpha}+(-1)^{\alpha+\gamma}U^{\alpha}_{~\gamma}T^{\gamma}_{~\alpha}(-1)^{\alpha}
=-T^{\alpha}_{~\delta}U^{\delta}_{~\alpha}(-1)^{\alpha}+(-1)^{\gamma}T^{\gamma}_{~\alpha}U^{\alpha}_{~\gamma}=0.\nn\eea
$\square$

%
%

\subsection{Trivialization (in)Dependence}\label{TiD}
Next, we need to check whether the definition of the Wilson line is independent of the trivialization of the graded bundle and whether it is independent of the exact position at which we choose to insert the transition function.

A change of trivialization of ${\cal E}$ is written as
\bea (x^A,\zeta_{\alpha})\to (x^A,\tilde \zeta_{\tilde{\alpha}})=(x^A,\zeta_{\beta}\Omega^{\beta}_{~\tilde{\alpha}}(x)).\nn\eea
Under this change
\bea T^{\tilde\beta}_{~\tilde\alpha}=(-1)^{\tilde\beta+\beta}\big[(\Omega^{-1})^{\tilde\beta}_{~\beta}Q\Omega^{\beta}_{~\tilde\alpha}\big]
+(-1)^{\tilde\beta+\alpha}(\Omega^{-1})^{\tilde\beta}_{~\alpha}T^{\alpha}_{~\beta}\Omega^{\beta}_{~\tilde\alpha}.\label{gauge_trans}\eea
Writing this equation slightly differently
\bea(-1)^{\tilde\gb}\Omega^{\gd}_{~\tilde\gb}T^{\tilde\beta}_{~\tilde\alpha}=(-1)^{\gd}Q\Omega^{\gd}_{~\tilde\alpha}
+(-1)^{\gd}T^{\gd}_{~\beta}\Omega^{\beta}_{~\tilde\alpha}\label{checked_later}\eea
shows that the gauge transformation rule is neatly written as
\bea \hat Q\Omega=0,~~~~\Omega=\zeta_{\gb} \Omega^{\gb}_{~\tilde\ga}\tilde\zeta^{\tilde\ga},\label{zero_leng_W}\eea
where $\tilde\zeta^{\tilde\ga}$ is the dual coordinate of $\tilde\zeta_{\tilde\ga}$ just as $\zeta^{\ga}$ is the dual of $\zeta^{\ga}$.

For an infinitesimal transformation $\Omega^{\alpha}_{~\tilde\beta}=\delta^{\alpha}_{\tilde\beta}+\epsilon^{\alpha}_{~\tilde\beta}$,
\bea \delta_{\epsilon}T^{\beta}_{~\alpha}=Q\epsilon^{\beta}_{~\alpha}
+T^{\beta}_{~\gamma}\epsilon^{\gamma}_{~\alpha}-(-1)^{\beta+\gamma}\epsilon^{\beta}_{~\gamma}T^{\gamma}_{~\alpha}.\label{gauge_transform}\eea
If one goes back to ordinary vector bundles and take $Q$ to be the de Rham differential, the formula above is just the gauge transformation of a connection $T$.

\begin{Prop}\label{Prop_triv_ind}
The change of the Wilson line under an infinitesimal change of trivialization is
\bea\delta_{\epsilon}U^{\alpha}_{~\gb}=\Big(-\epsilon^{\alpha}_{~\gd}(1)U^{\gd}_{~\gb}(1,0)+U^{\alpha}_{~\gd}(1,0)\epsilon^{\gd}_{~\gb}(0)\Big)
+(-1)^{\ga}\delta_BV^{\ga}_{~\gb}+(-1)^{\ga}T^{\ga}_{~\gd}(1)V^{\gd}_{~\gb}+(-1)^{\gd}V^{\ga}_{~\gd}T^{\gd}_{~\gb}(0),\label{triv_ind}\eea
where
\bea V^{\ga}_{~\gb}=\int_0^1 dt~U^{\alpha}_{~\gd}(1,t)\big[(-1)^{\gd}\partial_{\theta}
\BS{\epsilon}^{\gd}_{~\gamma}
\big]\Big|_{\theta=0}(t)U^{\gamma}_{~\gb}(t,0).\nn\eea
Or simply
\bea \delta_{\epsilon}U_C=\delta_B(\zeta_{\ga}V^{\ga}_{~\gb}\zeta^{\gb}),\label{triv_ind_simple}\eea
\end{Prop}
The proof is given in the appendix $\square$.

As mentioned earlier the Wilson-line maybe formed by combining two sections using a transition function across two trivialization patches. It suffices to look at the case of two patches
\bea U_C=\zeta_{\ga}U^{\ga}_{~\gc}(1,t)\cdot\Omega^{\gc}_{~\tilde\gd}(t)\cdot U^{\tilde \gd}_{~\tilde\gb}(t,0)\tilde\zeta^{\tilde \gb}\nn\eea
Under a change of trivialization in individual patches, the transition function obviously transforms covariantly, and using Eq.\ref{zero_leng_W},\ref{concise} we have
\bea \gd_{\ge}W_C=\gd_B\big(\zeta_{\ga}V^{\ga}_{~\gc}(1,t)\cdot\Omega^{\gc}_{~\tilde\gd}(t)\cdot U^{\tilde \gd}_{~\tilde\gb}(t,0)\tilde\zeta^{\tilde \gb}
+\zeta_{\ga}U^{\ga}_{~\gc}(1,t)\cdot\Omega^{\gc}_{~\tilde\gd}(t)\cdot V^{\tilde \gd}_{~\tilde\gb}(t,0)\tilde\zeta^{\tilde \gb}\big),\nn\eea
thus we see the insertion of transition function in between sections of Wilson-lines causes no further complication.

From the prop.\ref{Prop_triv_ind} and the above discussion we have an easy corollary
\begin{Cor}
If $\epsilon^{\ga}_{~\gb}$ is of degree zero, then the Wilson line is covariant under the change of trivialization.
\end{Cor}
This is because the term $\partial_{\theta}\BS{\epsilon}^{\gd}_{~\gamma}$ will have negative degree, and whether or not it is hit by the ensuing $\gd_B$ it lands in the ideal $(\textrm{eom},x^m_{p,<0}(y_0))$, see Eq.\ref{ideal}.

Take the setting from Ex.\ref{EX_Flat_bundle}, in particular Eq.\ref{Q_flat_bundle}, then the change of trivialization matrix $\epsilon^i_{~j}$ only depends on $M$, leading to
\bea \gd_{\ge}U(1,0)=-\ge(1)U(1,0)+U(1,0)\ge(0),\nn\eea
whose finite version is the familiar statement that the Wilson-line for a usual vector bundle transforms covariantly
\bea U(1,0)\to \Omega^{-1}(1)U(1,0)\Omega(0),\nn\eea
which is true even for non-flat connections.

Another important consequence of prop.\ref{Prop_triv_ind} is
\begin{Cor}The Wilson-loops are independent of the trivializations\end{Cor}
\noindent proof: From Eq,\ref{triv_ind} we get (remembering $x(0)=x(1)$, we write $T(0)=T(1)=T$ and $\ge(0)=\ge(1)=\ge$)
\bea \gd_{\ge}(-1)^{\ga}U^{\ga}_{~\ga}=-(-1)^{\ga}\epsilon^{\alpha}_{~\gd}U^{\gd}_{~\ga}+(-1)^{\ga}U^{\alpha}_{~\gd}\epsilon^{\gd}_{~\ga}
+(-1)^{\ga}\delta_BV^{\ga}_{~\ga}+T^{\ga}_{~\gd}V^{\gd}_{~\ga}+(-1)^{\ga+\gd}V^{\ga}_{~\gd}T^{\gd}_{~\ga},\nn\eea
the $\delta_BV^{\ga}_{~\ga}$ term drops since every $\ge^{\ga}_{~\gb}$ appearing in $V$ will be of zero degree. The rest of the terms also cancel after carefully renaming the dummy indices $\square$.
\begin{rmk}
This corollary shows that the Wilson-loops give bona fide functions on $\bar{\cal S}_{\textrm{O}}=\HomM(T[1]S^1,\M)$.
\end{rmk}

The trivialization dependence issue for the Wilson-lines must be dealt with before there is any sense of continuing, the general formula Eq.\ref{triv_ind} shows that one can either demand that $\ge^{\ga}_{~\gb}$ be of deg 0. This condition can be achieved by choosing a connection for $\E$, for example, in the case of Ex.\ref{Adjoint_rep}, we can choose a connection for $L\to M$, and split $T^*[2]L[1]$ non-canonically into $T^*[2]M\oplus L^*[1]\oplus L[1]$, then any further transition function will only depend on $M$ rather than the whole $L[1]$. In general, if we denote by $|\M|$ the degree zero part (the body) of $\M$, then it is always possible to split $\E\to\M$ into bundles over $|M|$ pulled back to $\M$. But we deem this approach \emph{unnatural}, since if one looks at the representation theory of $Q$-structure from the point of view of $Q$-equivariant graded vector bundle as in sec.\ref{RoQS}, one should be given the full freedom of choosing any trivialization. Besides, even if one can force $\ge$ to be degree zero by means of a connection, one must still confront the question of 'connection independence', which leads one to square one again. Thus we reject this solution.
{\color{black}\begin{rmk} The statements in ref.\cite{AbadSchaetz} prop.4.4 and cor.4.6 about trivialization independence seemed to us somewhat imprecise. Cor.4.6 states (in the nomenclature of this paper) that the Wilson-line transforms covariantly provided the iterated integral is such that it gives zero if any of the entry in the bar complex is \emph{invertible}. And they used prop.3.26 (which corresponds to rmk.\ref{rmk_neg_deg} of this paper) to justify this property. The iterated integral is such that it gives zero if any of the entry is of \emph{degree 0} (property 2 in prop.3.26), yet, $\Omega$ above, though invertible, may be of positive degree since it is a matrix. So their statement is weakened to the scenario of the first option above, namely covariance is only obtained when $\Omega^{\ga}_{~\gb}$ depends on the reduced manifold $|\M|$ alone. Consequently their discussion cannot be used to tackle the connection independence problem. \end{rmk}}

One may also be tempted to consider the $\gd_B$-cohomology since all the transgressing terms in Eq.\ref{triv_ind_simple} seems $\gd_B$ exact, but unfortunately the term $\zeta_{\ga}V^{\ga}_{~\gb}\zeta^{\gb}$ hit by $\gd_B$ is also trivialization dependent, thus not really $\gd_B$-exact. Having spurned the first two options, the only one left is, though unpalatable, that we must restrict the type of test manifolds $\Pa$ to have zero $Q$-structure. To see this, notice that by using the functor of points perspective, we are probing
the bundle $\E\textrm{nd}$ over $\bar{\cal S}$ (see the map in fig.\ref{ideal_situation}) with the test manifold $\Pa$, and naturally $\gd_B$ is mapped to $d^y$. By using only $\Pa$ with zero $Q$ (namely $N$-manifolds) as probes, we become totally oblivious of the $Q$-structure of $\E\textrm{nd}$. Simply put, rhs of Eq.\ref{triv_ind_simple} will eventually turn into $d^y$ and if $d^y=0$, we have that the Wilson line $U$ transforms covariantly. Furthermore, when $d^y=0$, we have
\bea \gd_B\BS{x}^m(t,\theta,y)=\BS{Q}^m(t,\theta,y)-D\BS{x}^m(t,\theta,y)=0,\nn\eea
in all degrees.

Note that by 'forgetting the $Q$-structure' we do not mean the forgetful functor that relegates one from the category of $NQ$ manifolds to $N$ manifolds. Rather, we are still working in the category of $NQ$-Mfld's. The next example shows us how much lies outside the detection range of an $N$-manifold probe.
\begin{Exa}\label{information_loss} Taking the setting from ex.\ref{Courant}, then
probing $\M=T[1]M$ with an $N$-manifold $\Pa$ will give us full information about $M$ but nothing about the fibre generator $v^{\mu}$, while a map $\Pa\to T^*[2]T[1]M$ that covers the map $\Pa\to T[1]M$ will know about the fibre generator $q_{\mu}$ but is unaware of the generator $p_{\mu}$, since both $v^{\mu}$ and $p_{\mu}$ are in the image of the $\hat D$ of Eq.\ref{D_hat}.

In general, probing $\M=L[1]$ with an $N$-manifold will only detect $M$ and those sections of $L$ in the kernel of the anchor. And we can crudely say that the Wilson-line is a morphism between the right kernel of $T(0)$ to the left kernel of $T(1)$.
\end{Exa}
Thus the information loss is quite severe; perhaps a better solution is called for to circumvent the trivialization dependence problem.

\subsection{Homotopy (in)Variance}\label{sec_HiV}
We next study how does the Wilson line change under a deformation of the super curve. By deformation, we mean those of type Eq.\ref{homotopy_abstract}, which we write again for convenience
(do pay attention that the \dga~$\textsf{B}$ in the setting \ref{homotopy_setting} is now generated by $t,\theta$ \emph{and} $y$)
\bea
&&\partial_s \BS{x}^A(t,\theta,y,s)=(\bar{\BS{x}}^C\frac{\partial}{\partial x^C}\gd\BS{x}^A)(t,\theta,y,s)+(D+d^y)\bar{\BS{x}}^A(t,\theta,y,s),\label{again}\eea
We shall write $\iota_{\bar x}$ for $\bar{\BS{x}}^C\partial_{x^C}$. The homotopy should also be extended to act on the generators of the section $\Gamma(\E)$ at $x(1)$ and $x(0)$, namely $\zeta_{\ga}$ and $\zeta^{\ga}$
\bea \partial_s\zeta_{\ga}=\zeta_{\gb}(\iota_{\bar x}T^{\gb}_{~\ga}(1)),~~~~\partial_s\zeta^{\ga}=(\iota_{\bar x}T^{\ga}_{~\gb}(0))\zeta^{\gb},\label{deform_fibre}\eea
straightforwardly from Eq.\ref{again}. As a side comment here, without the functor of points perspective, we find it quite hard to formulate the concept of deformation.

Thus quite clearly we have
\begin{Prop}\label{general_homotopy}
\bea \partial_sU^{\ga}_{~\gb}=d^y\iota_{\bar x}U^{\ga}_{~\gb}+\iota_{\bar x}\Big(-T^{\ga}_{~\gd}(1)U^{\gd}_{~\gb}(1,0)+(-1)^{\ga+\gc}U^{\ga}_{~\gc}(1,0)T^{\gc}_{~\gb}(0)\Big)\label{deform_not_simple}\eea
or simply
\bea\partial_sU_C=d^y\big((-1)^{\ga}\zeta_{\ga}(\iota_{\bar x}U^{\ga}_{~\gb})\zeta^{\gb}\big).\label{deform_simple}\eea
\end{Prop}
Proof: Using Eq.\ref{again} we have
\bea \partial_sU^{\ga}_{~\gb}(1,0)&=&\int_0^1dt~U^{\ga}_{~\gc}(1,t)\Big(-\partial_{\theta}\big[(\iota_{\bar x}\gd_B\BS{x}^A+d^y\bar {\BS{x}}^A)\partial_A\BS{T}\big]\Big|_{\theta=0}\Big)^{\gc}_{~\gd}(t)U^{\gd}_{~\gb}(t,0)\nn\\
&=&\{\gd_B,\iota_{\bar x}\}U^{\ga}_{~\gb}(1,0)+\int_0^1dt~U^{\ga}_{~\gc}(1,t)\Big(-\partial_{\theta}\big[(d^y\bar {\BS{x}}^A)\partial_A\BS{T}\big]\Big|_{\theta=0}\Big)^{\gc}_{~\gd}(t)U^{\gd}_{~\gb}(t,0).\nn\eea
Using $(\gd_B-d^y)\BS{x}^A=0$, we get 
\bea \partial_sU^{\ga}_{~\gb}(1,0)&=&d^y\iota_{\bar x}U^{\ga}_{~\gb}(1,0)+\iota_{\bar x}\gd_BU^{\ga}_{~\gb}(1,0)\nn\\
&=&d^y\iota_{\bar x}U^{\ga}_{~\gb}(1,0)+\iota_{\bar x}\Big(-T^{\ga}_{\gd}(1)U^{\gd}_{~\gb}(1,0)+(-1)^{\ga+\gc}U^{\ga}_{~\gc}(1,0)T^{\gc}_{~\gb}(0)\Big)\nn\eea
\begin{flushright}
$\square$
\end{flushright}
We can illustrate the formula given above to the simple example of Courant algebroid, indeed the simplest examples are often the trickiest to get straight, but we leave them to the appendix to curb the length of this section.

The reparameterization is also a type of homotopy, we have the following important
\begin{cor} The Wilson-line is invariant under re-parametrization of $[0,1]$\end{cor}\label{cor_repara}

Proof: Take $\epsilon(t)\geq0$ and $\ge(0)=\ge(1)=0$ be the vector field generating the re-parametrization, then
we can choose $\bar{\BS{x}}=\ge(t)\partial_{\theta}\BS{x}$, for we have from Eq.\ref{again}
\bea \partial_s \BS{x}(t,\theta,y,s)&=&\ge(t)\partial_{\theta}\big(\gd\BS{x}(t,\theta,y,s)\big)+(D+d^y)\ge(t)\partial_{\theta}{\BS{x}}(t,\theta,y,s),\nn\\
&=&\ge(t)\partial_{\theta}\Big(\gd\BS{x}(t,\theta,y,s)-(D+d^y){\BS{x}}(t,\theta,y,s)\Big)+\Big(\ge(t)\partial_{t}+\dot{\ge}(t)\theta\partial_{\theta}\Big)
{\BS{x}}(t,\theta,y,s)\nn\\
&=&\Big(\ge(t)\partial_{t}+\dot{\ge}(t)\theta\partial_{\theta}\Big){\BS{x}}(t,\theta,y,s)=L_{\ge}{\BS{x}}(t,\theta,y,s),\nn\eea
where $L_{\ge}$ is the Lie derivative on [0,1]. Finally the conclusion follows from $\partial_{\theta}^2=0$ $\square$

To complete the discussion of invariance under homotopy, we need to consider a Wilson-line that traverses a number of trivialization patches. The problem is quite simple, consider the homotopy of the map $*\times\Pa\to\M$, where $*$ is a connecting point of two pieces of Wilson-lines. It is clear then under the homotopy
\bea\partial_s\Omega^{\ga}_{~\tilde\gb}(x)&=&\iota_{\bar x}(Q\Omega^{\ga}_{~\tilde\gb}(x))+d^y\iota_{\bar x}\Omega^{\ga}_{~\tilde\gb}(x),\nn\\
\textrm{equivalently}~~~\partial_s\Omega(x)&=&d^y\iota_{\bar x}\Omega(x).\label{deform_trans}\eea
where $\Omega(x)$ is as defined in Eq.\ref{zero_leng_W} $\Omega(x)=\zeta_{\ga}\Omega^{\ga}_{~\tilde\gb}(x)\tilde\zeta^{\tilde\gb}$ and we have used $\hat Q\Omega=0$. The last equation has the same form as Eq.\ref{deform_simple}, actually one may well think of the transition function as a \emph{zero-length Wilson-line}, as such, the spliced Wilson line will behave the same way under homotopy according to Eq.\ref{deform_simple}.

Having sorted out this issue, we can derive the crucial property of the Wilson-loops from the general formula prop.\ref{general_homotopy}
\begin{cor}\label{Prop_Deform0}
Wilson loop is invariant under the homotopy of the map $T[1]S^1\times\Pa\to\M$
\end{cor}
\noindent proof: Remembering $x(1)=x(0)$, we have
\bea \partial_sU_{\textrm{O}}=(-1)^{\ga}\partial_sU^{\ga}_{~\ga}=(-1)^{\ga}d^y\iota_{\bar x}U^{\ga}_{~\ga}+(-1)^{\ga}\iota_{\bar x}\Big(-T^{\ga}_{\gd}U^{\gd}_{~\ga}+(-1)^{\ga+\gc}U^{\ga}_{~\gc}T^{\gc}_{~\ga}\Big)\nn\eea
The first term on rhs drops since $\iota_{\bar x}U^{\ga}_{~\ga}$ has degree $-1$. The rest also cancel
\bea -(-1)^{\ga}T^{\ga}_{\gd}U^{\gd}_{~\ga}+(-1)^{\gc}U^{\ga}_{~\gc}T^{\gc}_{~\ga}=-(-1)^{\ga}T^{\ga}_{\gd}U^{\gd}_{~\ga}+(-1)^{\gc}T^{\gc}_{~\ga}U^{\ga}_{~\gc}
=-(-1)^{\ga}T^{\ga}_{\gd}U^{\gd}_{~\ga}+(-1)^{\ga}T^{\ga}_{~\gc}U^{\gc}_{~\ga}=0.\nn\eea
\begin{rmk}
Cor.\ref{Prop_Deform0} shows that the Wilson-loops are very special functions on $\bar{\cal S}$, in that they are invariant under homotopies. Thus if the quotient of $\bar{\cal S}$ by homotopy makes sense, then the Wilson-loops will go down to the quotient.\end{rmk}


%
%


%
It pays to look at an example that shows roughly what is happening to the Wilson-line under deformation, yet unincumbered by the petty signs.
\begin{Exa}\label{Wilson_flat_bundle}
Take the setting from Ex.\ref{EX_Flat_bundle}, and let $\phi(t)$ denote the mapping $[0,1]\to M$, and $H(t,s)$ a homotopy of $\phi(t)$ with $H(t,0)=\phi(t)$. And we denote the infinitesimal deformation by a vector field $v(t)=\partial_sH(t,s)\big|_{s=0}$, vanishing at the two ends $v(1)=v(0)=0$. Thus the change of $U(0,1)$ due to this deformation is
\bea \partial_sU(1,0)=\partial_s\BB{P}\exp\Big\{\int_0^1 (H^*A)(t,s)\Big\}\Big|_{s=0}=\int_0^1dt~U(1,t)\big((\phi^*L_vA)(t)\big)~U(t,0)\nn\eea
where $A$ is the flat connection.
\bea &&L_vA=\iota_vdA+d\iota_vA=-[\iota_vA,A]+d\iota_vA\nn\\
\Rightarrow&& \phi^*(L_vA)=\phi^*(-[\iota_vA,A])+d\phi^*\iota_vA,\nn\eea
when confusion is unlikely to arise, we will drop $\phi^*$.
\bea \partial_s U(1,0)=\int dt~U(1,t)\Big(-[\iota_vA,A]+\frac{d}{dt}\iota_vA\Big)(t)~U(t,0)\nn\eea
Integrate by part the second term and use the relation
\bea \frac{d}{dt_1}U(t_1,t_0)=A(t_1)U(t_1,t_0),~~\frac{d}{dt_0}U(t_1,t_0)=-U(t_0,t_1)A(t_0)\nn\eea
We get
\bea\partial_s U(1,0)=\int dt~U(1,t)\Big(-[\iota_vA,A]-A\iota_vA+(\iota_vA)A\Big)(t)~U(t,0)=0,\nn\eea
which is the expected result when deforming the Wilson line for a flat bundle.
\end{Exa}

{\color{black}The final proposition summarizes the discussion of this section.
\begin{Prop}
For a representation of an $NQ$-mfld $\M$, which is a $Q$-equivariant bundle $\E\to\M$, one can form the Wilson line operator for a super-curve, given in Eq.\ref{wilson_line}. It is covariant under a change of trivialization if one probes it (in the sense of functor of points perspective) with a test manifold $\L$ with zero $Q$-structure. In such a case, the Wilson line is also invariant under homotopy deformation of the super-curve.

In contrast, if one takes a trace as in Eq.\ref{Wilson_loop}, the resulting Wilson-loop operator is a well-defined function on the presheaf $\HomM_{NQ}(T[1]S^1,\M)$, and furthermore the Wilson-loop is invariant under homotopy deformation of the super loop.\end{Prop}}

\section{Summary and Possible Applications}\label{sec_ingetnamn}

The central application of the lengthy discussion of the previous section about the invariance property of Wilson-lines/loops is the following.
Consider the mapping problem $\Hom_{NQ}(T[1]X,\M)$, where $X$ is a usual manifold; the mapping is given by the generalization of flat connections, as sketched in sec.\ref{MCSaI}. We study the class of functions of these flat connections spanned by the Wilson loops. We embed a super-loop $T[1]S^1$ into $T[1]X$ and denote by $LX$ the loop space of $X$ (formerly denoted as $\bar{\cal S}_{\textrm{O}}$)
\bea T[1]S^1\longrightarrow T[1]X \stackrel{\varphi}{\longrightarrow} \M,\label{physics}\eea
and the map $\varphi$ is given by the generalized flat connections.

By picking a representation, a $Q$-equivariant graded bundle over $\M$, we can form the Wilson-loop which is a function that depends on $LX$ and also on the mapping $\varphi:\,\Hom_{NQ}(T[1]X,\M)$. By the earlier discussions, the Wilson-loop is independent of the trivialization thus it is a function of gauge equivalence class of generalized flat connections. The Wilson loop also only depends on the homotopy class of loops $S^1\to X$. This is nothing but the generalization of the familiar statement: gauge equivalence classes of flat connections are determined by holonomies.

In physics, this scenario Eq.\ref{physics} occurs often. We have a TFT whose fields are the mapping Eq.\ref{physics}, taking $X=\Sigma\times\BB{R}$ one can go to the Hamiltonian formalism. The reduced phase space are the gauge equivalence class of maps $T[1]\Sigma\to\M$. Take the Chern-Simons theory as an example, for which $\M=\FR{g}[1]$, and take $X=\Sigma_2\times \BB{R}$, the reduced phase space the moduli space of flat connections on the Riemann surface $\Sigma_2$. We have seen that the Wilson loops are functions on this moduli space. In the works \cite{Andersen} \cite{FockRosly}, it was shown that for Riemann surfaces with punctures and certain gauge groups, the Wilson loops corresponds to  complete basis of algebraic functions of the moduli space of flat connections, thus one can study directly the algebraic structure of Wilson-loops to study the moduli space. The proof of completeness uses the Peter-Weyl theorem combined with the classic invariant theory, this route will not be viable for us to take, since we only have but a scanty grip of the possible representations of an $NQ$-manifold. And the results of prop.\ref{iso_bar_fixed}, \ref{iso_bar_loop} are but a far cry toward answering whether the Wilson-loops in an $NQ$-manifold gives a complete basis or not. Thus we have plenty of room for future work here.  One can also look at the Poisson $\sigma$-model, then our approach resembles that of Cattaneo and Felder \cite{2000math3023C}, who tackled the problem of integration of Poisson structure by looking at the reduced phase space of $T[1][0,1]\to \M$, where $\M$ is $T^*[1]M$ equipped with a $Q$-structure constructed from the Poisson structure of $M$.

Compared to the Wilson-loops which are well-defined objects, the Wilson lines do not transform covariantly under a change of trivializations for arbitrary test manifolds. The Wilson-lines do enjoy the invariance under re-parametrizations, thus if we pick
two super curves $C_1,\,C_2$ with $C_1(1)=C_2(0)$, then
\bea U_{C_2}U_{C_1}=U_{C_2\circ C_1}.\label{composition}\eea
We would like to point out a few analogies with the familiar problem of integration of a Lie algebra $\int\FR{g}[1]$. Let the function $T^{\ga}_{~\gb}(g)$ be a matrix representation of $G$, i.e.\footnote{Suppose the group $G$ integrating $\FR{g}$ is compact, then it is known by Peter-Weyl theorem that the set of matrix coefficients $T^{\ga}_{~\gb}(g)$ for finite dimensional representations is dense in the space of continuous complex functions on $G$}
\bea T^{\ga}_{~\gc}(g_1)T^{\gc}_{~\gb}(g_2)=T^{\ga}_{~\gb}(g_1g_2),\nn\eea
also consider a map $M\times M\stackrel{\varphi}{\to} G$ satisfying the condition Eq.\ref{G_descent}, then the Wilson line $U^{\ga}_{~\gb}$ associated to a curve starting and ending at $x,y\in M$ is the pull back of the function $T^{\ga}_{~\gb}$ on $G$ to $M\times M$ under the map $\varphi$. In the graded setting, the integration object
$\int \M$ has thus far been understood as a presheaf Eq.\ref{significant} and beyond this we have little knowledge of it. However, as the motto of graded geometry, one should study $\int\M$ by looking at the set morphisms from other graded manifolds to it. In particular, we can construct maps
\bea \So\times \So\to \smalint\M,\nn\eea
with the same property as Eq.\ref{G_descent}. Concretely, for an $NQ$-manifold $\M$ with representation $\E\to\M$, one first solves the mapping problem $\Hom_{NQ}(T[1]\So,\M)$, which then gives the set of generalized flat connections. For two points in $\So$, pick a curve connecting the two and construct the Wilson-line operator associated to this curve. Due to its invariance under gauge transformation, the Wilson-line only depends on the initial and final point and the composition property Eq.\ref{composition} will guarantee Eq.\ref{G_descent}. We can interpret similarly that our graded Wilson-lines are certain functions over $\int\M$ pulled back to $\So\times \So$ as in the group case.

In the work of Abad and Crainic \cite{AbadCrainic}, the representation of groupoid up to homotopy was defined, for which the composition law fails by a $T_0$-exact term (see Eq.\ref{T-square} for the notation). But in our setup $T_0$-exact terms will be set to zero since we already required $d^y=0$. It is perhaps possible to relax the condition on $d^y$ somehow and recover Abad and Crainic's representation of groupoid up to homotopy from our formalism.  Finally let us point out the relationship between our work and that of ref.\cite{AbadSchaetz}. The authors there used a sophisticated parameterization due to Igusa of paths from the first vertex to the last vertex of a simplex $\Delta^n$ fitted inside $\M$ (in the sense of $NQ$-manifold morphism of course) by a cube $[0,1]^{n-1}$. The component of Wilson-line of degree $n-1$ can be integrated over the cube. Thus they obtained a representation of the simplicial set $(\int \Ta)_n$, which they call the representation up to homotopy of the $\Pi_{\infty}$-groupoid. In comparison, we did not integrate over $[0,1]^{n-1}$ and retained the components of Wilson-lines of all degrees as functions on $\bar{\cal S}$, modulo the trivialization issues. To further compare with their work, such as studying the quasi-isomorphisms of two representations, one needs to fix a connection which will not be taken up in this paper.

\bigskip\bigskip

\noindent{\bf\Large Acknowledgement}:
\bigskip

{\color{black}The authors would like to thank E.Getzler for discussion on the problem of trivialization dependence.}
M.Z. thanks INFN Sezione di Firenze, Universit\`a
di Firenze and  KITP, Santa Barbara where part of this work was carried out.
The research of M.Z. is supported by VR-grant 621-2008-4273  and
 was supported in part by DARPA under Grant No.
HR0011-09-1-0015 and by the National Science Foundation under Grant
No. PHY05-51164.

\bigskip\bigskip

\appendix
\section{Representation up to Homotopy in Graded Language}
For a Lie algebroid data given as in Ex.\ref{Lie_algebroid}, one builds a \emph{representation up to homotopy} ($\rep$) \cite{abad-2009} on a complex of vector bundles $E_{\sbullet}$ over $M$. The representation is given by a differential $D$ with the following decomposition
\bea D=D_{0,1}+D_{1,0}+D_{2,-1}+\cdots\nn\eea
similar to the decomposition Eq.\ref{decomp_T}.
If we denote $\Omega^{p,q}(L,E)=\Gamma(\wedge^pL^*\otimes E_q)$, then the effect of $D$ is
\bea D_{s,1-s}:~\Omega^{p,q}(L,E)\to\Omega^{p+s,q+1-s}(L,E),~~~~s\geq0.\nn\eea

In the graded language, we represent the \emph{sections} of the complex $\Omega^{p,q}(L,E)$ as \emph{coordinates} of a graded manifold. In this way, we can simply formulate $D$ as a homological vector field $Q$ on this GM. For a Lie algebroid problem above, we first take
$\M=L[1]$, with a $Q$-structure defined as in Eq.\ref{Q_Lie_algebroid}. Thus the degree 1 fibre coordinate $\ell^A$ takes the place of $\Omega^{1,0}(L,E)$. And we pick further a $Q$-equivariant bundle $\E\stackrel{\pi}{\to} \M$ as in fig.\ref{def_rep}.
The transition function of the bundle $\E$ can be made to depend on $M$ if we pick a connection for $L$. Assuming this done, and that
$\E$ is split into bundles over $M$: $\E\sim \E^*\oplus L[1]$. Again, the degree $q$ fibre coordinate of $E^*$ takes the place of $\Omega^{0,q}(L,E)$, while in general a section of $\E$ gives the full $\Omega^{p,q}(L,E)$. With this said, the differential $D$ above is given by $\hat{Q}$
\bea &&D=\hat{Q},\nn\\
&&D_{1,0}=Q+T_0,~~D_{0,1}=T_1,~~ D_{2,-1}=T_2,\cdots\nn\eea
where the decomposition of $T$ was given in Eq.\ref{decomp_T}.

Next we see how this works concretely for the adjoint representation of a Lie algebroid (notation as in Ex.\ref{Adjoint_rep}).
Choose a local basis $e_A$ for the sections of $L$, we write the connection for $L$ in this basis as
\bea \nabla_{\mu}e_A=\partial_{\mu}e_A+\Gamma_{\mu A}^Be_B,~~~~e_A\in\Gamma(L).\nn\eea
Take $\E=T^*[2]L[1]$, and it can be split into
\bea T^*[2]L[1]\sim T^*[2]M\oplus L^*[1] \oplus L[1]\nn\eea
This is done by shifting the coordinate
\bea \tilde{p}_{\mu}=p_{\mu}+\Gamma^A_{\mu
B}\bar\ell_A\ell^B.\nn\eea
It is easy to check that $\tilde{p}_{\mu}$ transforms as $\partial_{\mu}$ by using the definition of $\Gamma_{\mu}$, and it is also clear that $\bar\ell_A,\ell^A$ transforms as sections of $L,L^*$. Thus we read off the complex of vector bundles as
\bea E_{\sbullet}=\cdots 0\to L\to TM\to 0\to\cdots,\nn\eea
with only two levels. In this case, the coordinate $\tilde p_{\mu}$ of $T^*[2]$ represents the section of $TM$ at level 2 while the coordinate $\bar\ell_A$ of $L^*[1]$ the section of $L$ at level 1.

Define some new shifted quantities
\bea \Gamma^A_{CB}=A_C^{\mu}\Gamma_{\mu B}^A,~~~~~\tilde f^A_{BC}=f^A_{BC}+\Gamma_{[CB]}^A\nn\eea
We observe that
the shifted structure function $\tilde
f^A_{BC}$ transforms homogeneously for all
three indices under a change of trivialization of $L$ hence it is a tensor.

To construct $\hat Q$, we note that $Q$ has a straightforward Hamiltonian lift
\bea \Theta=2p_{\mu}
A^{\mu}_C\ell^C-f^A_{BC}\bar\ell_A\ell^B\ell^C=2\tilde p_{\mu}A^{\mu}_C\ell^C-\tilde{f}^A_{BC}\bar\ell_A\ell^B\ell^C,\label{encapsulate}\eea
in the sense that $\{\Theta,f(x,\ell)\}=Qf(x,\ell)$. And $\hat Q$ is defined as $\hat Qf(x,\ell,\bar\ell,p)=\{\Theta,f(x,\ell,\bar\ell,p)\}$ which lifts $Q$ from $\M$ to $\E$. We see the point of using graded manifold language is that all is now encapsulated in one equation
\bea \hat{Q}^2=0\Leftrightarrow \{\Theta,\Theta\}=0,\nn\eea

Writing $\hat{Q}$ in terms of the shifted quantities we have
\bea \hat{Q}=Q+\Big(2\tilde p_{\mu}A_A^{\mu}-2\tilde f^C_{BA}\bar\ell_C\ell^B\Big)\frac{\partial}{\partial\bar\ell_A}
-\Big(2(\partial_{\mu}A_A^{\nu}-\Gamma_{\mu
A}^BA_B^{\nu})\tilde p_{\nu}\ell^A+(\nabla_{\mu}\tilde f_{AB}^C)\bar\ell_C\ell^A\ell^B\Big)\frac{\partial}{\partial\tilde p_{\mu}}.\nn\eea
Decompose
$\hat Q$ according to degrees
\bea D_{1,0}&=&Q-2\tilde f^C_{BA}\bar\ell_C\ell^B\frac{\partial}{\partial\bar\ell_A}-2(\partial_{\mu}A_A^{\nu}-\Gamma_{\mu
A}^BA_B^{\nu})\tilde p_{\nu}\ell^A\frac{\partial}{\partial\tilde p_{\mu}}\nn\\
D_{0,1}&=&\Big(2\tilde p_{\mu}A_A^{\mu}\Big)\frac{\partial}{\partial\bar\ell_A}\nn\\
D_{2,-1}&=&-(\nabla_{\mu}\tilde f_{AB}^C)\bar\ell_C\ell^A\ell^B\frac{\partial}{\partial\tilde p_{\mu}}.\label{D_decomp}\eea
In ref.\cite{abad-2009}, $D_{1,0},D_{0,1},D_{2,-1}$ are given names $\rho,\nabla^{bas},R^{bas}$ ($\rho$ is the anchor $L\to TM$) , and are constructed in coordinate-free manner as
\bea
&&\nabla_{\alpha}^{bas}(X)=\rho(\nabla_X(\alpha))+[\rho(\alpha),X]\nn\\
&&\nabla_{\alpha}^{bas}(\beta)=\nabla_{\rho(\beta)}(\alpha)+[\alpha,\beta],~~~~\alpha,\beta\in \Gamma(L),~~X\in\Gamma(T_M).\nn\eea
And $R^{bas}\in
\Gamma(\wedge^2 L^*\otimes \hom(T_M,L))$ is the exotic curvature
\bea
R^{bas}(\alpha,\beta)(X)=\nabla_X([\alpha,\beta])-[\nabla_X(\alpha),\beta]-[\alpha,\nabla_X(\beta)]
-\nabla_{\nabla^{bas}_{\beta}X}(\alpha)+\nabla_{\nabla^{bas}_{\alpha}X}(\beta).\nn\eea
To compare these expression to the ones given in local coordinates, we denote
\bea &&\nabla_{A\nu}^{\mu}=
\nabla^{bas}_{\bar\ell_A}(\partial_{\nu})\circ x^{\mu},~~~\nabla_{AB}^{C}=
\langle\nabla^{bas}_{\bar\ell_A}(\bar\ell_B),\ell^C\rangle,~~~R_{AB\mu}^{C}=\langle
R(\bar\ell_A,\bar\ell_B)(\partial_{\mu}),\ell^C\rangle\nn\eea
If one is so disposed, he can find the local expressions for these quantities by comparing them with Eq.\ref{D_decomp}, yet it is way more efficient to work with Eq.\ref{encapsulate} instead.

\section{Proof of some Propositions}
The Wilson line is defined as
\bea U^{\alpha}_{~\beta}(1,0)=\left(\BB{P}\exp\Big(-\int_0^1dtd\theta~\BS{T}\Big)\right)^{\alpha}_{~\beta}\nn\eea
associated with a super line $T[1][0,1]\to \M$. Again, the notation is that the boldface symbols are promoted to being superfields
$\BS{T}=T(\BS{x}(t,\theta))$.

\subsection{Proof of Prop.\ref{FIP}}
The BRST differential $\delta_B\BS{x}^A=-D\BS{x}^A+\BS{Q}^A$ acting on $U(1,0)$ gives
\bea \delta_B U^{\alpha}_{~\beta}(1,0)=\int_0^1dt~(-1)^{\alpha+\gamma}U^{\alpha}_{~\gamma}(1,t)
\Big(\partial_{\theta}\big[-D\BS{T}^{\gamma}_{~\delta}+\BS{QT}^{\gamma}_{~\delta}\big]\Big|_{\theta=0}\Big)(t)U^{\delta}_{~\beta}(t,0).\nn\eea

The main point is to integrate by part the first term in the brace, giving some surface term of type $\BS{T}\BS{T}$ which cancels the $\BS{QT}$ term in the same brace.
\bea
\delta_B U^{\alpha}_{~\beta}(1,0)=\int_0^1dt~(-1)^{\alpha+\gamma}U^{\alpha}_{~\gamma}(1,t)
\Big(-\partial_tT^{\gamma}_{~\delta}+\partial_{\theta}\big[\BS{QT}^{\gamma}_{~\delta}\big]\Big|_{\theta=0}\Big)(t)U^{\delta}_{~\beta}(t,0),\nn\eea
integrate by part and use the relations
\bea \partial_tU^{\alpha}_{~\beta}(t,0)=-\big(\partial_{\theta}\BS{T}\big|_{\theta=0}\big)^{\alpha}_{~\gamma}(t)U^{\gamma}_{~\beta}(t,0),~~~
\partial_tU^{\alpha}_{~\beta}(1,t)=U^{\alpha}_{~\gamma}(1,t)\big(\partial_{\theta}\BS{T}\big|_{\theta=0}\big)^{\gamma}_{~\beta}(t),\nn\eea
we get
\bea\delta_B U^{\alpha}_{~\beta}(1,0)&=&-T^{\alpha}_{~\delta}(1)U^{\delta}_{~\beta}(1,0)+(-1)^{\alpha+\gamma}U^{\alpha}_{~\gamma}(1,0)T^{\gamma}_{~\beta}(0)\nn\\
&&+\int_0^1dt~(-1)^{\alpha+\gamma}U^{\alpha}_{~\delta}(1,t)\Big(\big(\partial_{\theta}\BS{T}\big|_{\theta=0}\big)^{\delta}_{~\gamma}
T^{\gamma}_{~\epsilon}\Big)(t)~U^{\epsilon}_{~\beta}(t,0)\nn\\
&&-\int_0^1dt~(-1)^{\alpha+\gamma}U^{\alpha}_{~\gamma}(1,t)
\Big(T^{\gamma}_{~\delta}\big(\partial_{\theta}\BS{T}\big|_{\theta=0}\big)^{\delta}_{~\epsilon}\Big)(t)~U^{\epsilon}_{~\beta}(t,0)\nn\\
&&+\int_0^1dt~(-1)^{\alpha+\gamma}U^{\alpha}_{~\gamma}(1,t)\Big(\big(\partial_{\theta}\BS{QT}\big|_{\theta=0}\big)^{\gamma}_{~\delta}\Big)(t)
~U^{\delta}_{~\beta}(t,0)\nn\eea
The last three terms can be combined into
\bea
&&\int_0^1dt~U^{\alpha}_{~\delta}(1,t)\Big((-1)^{\alpha+\gamma}\big(\partial_{\theta}\BS{T}\big)^{\delta}_{~\gamma}
T^{\gamma}_{~\epsilon}-(-1)^{\alpha+\delta}T^{\delta}_{~\gamma}\big(\partial_{\theta}\BS{T}\big)^{\gamma}_{~\epsilon}
+(-1)^{\alpha+\delta}\big(\partial_{\theta}\BS{QT}\big)^{\delta}_{~\epsilon}
\Big)\big|_{\theta=0}(t)~U^{\epsilon}_{~\beta}(t,0)\nn\\
&=&(-1)^{\alpha}\int_0^1dt~U^{\alpha}_{~\delta}(1,t)\Big((-1)^{\gamma}\big(\partial_{\theta}\BS{T}^{\delta}_{~\gamma}\BS{T}^{\gamma}_{~\epsilon}\big)
+(-1)^{\delta}\big(\partial_{\theta}\BS{QT}\big)^{\delta}_{~\epsilon}
\Big)\big|_{\theta=0}(t)~U^{\epsilon}_{~\beta}(t,0)\nn\\
&=&(-1)^{\alpha}\int_0^1dt~U^{\alpha}_{~\delta}(1,t)\partial_{\theta}\Big((-1)^{\gamma}\BS{T}^{\delta}_{~\gamma}\BS{T}^{\gamma}_{~\epsilon}
+(-1)^{\delta}\BS{QT}^{\delta}_{~\epsilon}
\Big)\Big|_{\theta=0}(t)~U^{\epsilon}_{~\beta}(t,0)=0,\nn\eea
where Eq.\ref{normalized_MC} is used in the last step. The remaining terms is the BRST variation of a Wilson line
\bea \delta_B U^{\alpha}_{~\beta}(1,0)=-T^{\alpha}_{~\delta}(1)U^{\delta}_{~\beta}(1,0)+(-1)^{\alpha+\gamma}U^{\alpha}_{~\gamma}(1,0)T^{\gamma}_{~\beta}(0)
\label{BRST_WL}\eea
The proof is complete.

\subsection{Proof of Prop.\ref{Prop_triv_ind}}
As a crucial proposition of the paper, we present the proof of it with more details.
Under an infinitesimal change of trivialization $V^{\alpha}_{~\tilde\beta}=\delta^{\alpha}_{\tilde\beta}+\epsilon^{\alpha}_{~\tilde\beta}$, the matrix $T$ changes by
\bea \delta_{\tiny{\textrm{chg of triv}}}T^{\beta}_{~\alpha}=Q\epsilon^{\beta}_{~\alpha}
+T^{\beta}_{~\gamma}\epsilon^{\gamma}_{~\alpha}-(-1)^{\beta+\gamma}\epsilon^{\beta}_{~\gamma}T^{\gamma}_{~\alpha}.\eea
The change of the Wilson-line is
\bea \delta_{\tiny{\textrm{chg of triv}}}U^{\alpha}_{~\rho}(1,0)=\int_0^1 dt~U^{\alpha}_{~\beta}(1,t)\Big\{-\partial_{\theta}
\big[\underbrace{\BS{Q\epsilon}^{\beta}_{~\gamma}}_a
+\underbrace{\BS{T}^{\beta}_{~\delta}\BS{\epsilon}^{\delta}_{~\gamma}}_{3}-\underbrace{(-1)^{\beta+\delta}\BS{\epsilon}^{\beta}_{~\delta}\BS{T}^{\delta}_{~\gamma}}_{6}\big]\Big|_{\theta=0}\Big\}(t)U^{\gamma}_{~\rho}(t,0).\nn\eea
Write $Q\epsilon$ in term-$a$ as $D\epsilon+\delta_{B}\epsilon$,
\bea a=\int_0^1 dt~U^{\alpha}_{~\beta}(1,t)\Big\{-\partial_{\theta}
\big[\underbrace{D\BS{\epsilon}^{\beta}_{~\gamma}}_b+\underbrace{\delta_{B}\BS{\epsilon}^{\beta}_{~\gamma}}_c
\big]\Big|_{\theta=0}\Big\}(t)U^{\gamma}_{~\rho}(t,0)\nn\eea
Integrate by part the term-$b$
\bea b&=&-\epsilon^{\alpha}_{~\beta}(1)U^{\beta}_{~\rho}(1,0)+U^{\alpha}_{~\beta}(1,0)\epsilon^{\beta}_{~\rho}(0)\nn\\
&&+\int_0^1 dt~U^{\alpha}_{~\delta}(1,t)\Big\{\underbrace{\big[\partial_{\theta}\BS{T}^{\delta}_{~\beta}
\big]\Big|_{\theta=0}\epsilon^{\beta}_{~\gamma}}_{1}\Big\}(t)U^{\gamma}_{~\rho}(t,0)
-\int_0^1 dt~U^{\alpha}_{~\beta}(1,t)\Big\{\underbrace{\epsilon^{\beta}_{~\gamma}\big[\partial_{\theta}\BS{T}^{\gamma}_{~\delta}
\big]\Big|_{\theta=0}}_{4}\Big\}(t)U^{\delta}_{~\rho}(t,0).\nn\eea
Pull out $\gd_B$ from term-$c$
\bea c&=&(-1)^{\alpha+\beta}\delta_B\int_0^1 dt~U^{\alpha}_{~\beta}(1,t)\Big\{\partial_{\theta}
\big[\BS{\epsilon}^{\beta}_{~\gamma}
\big]\Big|_{\theta=0}\Big\}(t)U^{\gamma}_{~\rho}(t,0)\nn\\
&&-(-1)^{\alpha+\beta}\int_0^1 dt~(\delta_BU^{\alpha}_{~\beta}(1,t))\Big\{\partial_{\theta}
\big[\BS{\epsilon}^{\beta}_{~\gamma}
\big]\Big|_{\theta=0}\Big\}(t)U^{\gamma}_{~\rho}(t,0)\nn\\
&&+(-1)^{\beta+\gamma}\int_0^1 dt~U^{\alpha}_{~\beta}(1,t)\Big\{\partial_{\theta}
\big[\BS{\epsilon}^{\beta}_{~\gamma}
\big]\Big|_{\theta=0}\Big\}(t)\delta_BU^{\gamma}_{~\rho}(t,0)\nn\\
&=&\delta_B(\cdots)-(-1)^{\alpha+\beta}\int_0^1 dt\Big\{-T^{\alpha}_{~\delta}(1)U^{\delta}_{~\beta}(1,t)+\underbrace{(-1)^{\alpha+\delta}U^{\alpha}_{~\delta}(1,t)T^{\delta}_{~\beta}}_{2}\Big\}\big[\partial_{\theta}
\BS{\epsilon}^{\beta}_{~\gamma}
\big]\Big|_{\theta=0}(t)U^{\gamma}_{~\rho}(t,0)\nn\\
&&+(-1)^{\beta+\gamma}\int_0^1 dt~U^{\alpha}_{~\beta}(1,t)\big[\partial_{\theta}
\BS{\epsilon}^{\beta}_{~\gamma}
\big]\Big|_{\theta=0}(t)\Big\{\underbrace{-T^{\gamma}_{~\delta}(t)U^{\delta}_{~\rho}(t,0)}_{5}
+(-1)^{\gamma+\delta}U^{\gamma}_{~\delta}(t,0)T^{\delta}_{~\rho}(0)\Big\}\nn\eea
Collect everything together we have
\bea
&&\delta_{\tiny{\textrm{chg of triv}}}U^{\alpha}_{~\rho}(1,0)\nn\\
&=&\int_0^1 dt~U^{\alpha}_{~\delta}(1,t)\Big\{\underbrace{(\partial_{\theta}\BS{T}^{\delta}_{~\beta})\epsilon^{\beta}_{~\gamma}}_{1}
-\underbrace{(-1)^{\beta+\delta}T^{\delta}_{~\beta}\partial_{\theta}\BS{\epsilon}^{\beta}_{~\gamma}}_{2}
-\underbrace{(\partial_{\theta}\BS{T}^{\delta}_{~\beta}\BS{\epsilon}^{\beta}_{~\gamma})}_{3}
\Big\}\Big|_{\theta=0}(t)U^{\gamma}_{~\rho}(t,0)\nn\\
&&+\int_0^1 dt~U^{\alpha}_{~\beta}(1,t)\Big\{\underbrace{-\epsilon^{\beta}_{~\gamma}\partial_{\theta}\BS{T}^{\gamma}_{~\delta}}_{4}
-\underbrace{(-1)^{\beta+\gamma}(\partial_{\theta}\BS{\epsilon}^{\beta}_{~\gamma})T^{\gamma}_{~\delta}}_{5}
+\underbrace{(-1)^{\beta+\gamma}\partial_{\theta}(\BS{\epsilon}^{\beta}_{~\gamma}\BS{T}^{\gamma}_{~\delta})}_{6}
\Big\}\Big|_{\theta=0}(t)U^{\delta}_{~\rho}(t,0)\nn\\
&&+\delta_B(\cdots)-\epsilon^{\alpha}_{~\beta}(1)U^{\beta}_{~\rho}(1,0)+U^{\alpha}_{~\beta}(1,0)\epsilon^{\beta}_{~\rho}(0)\nn\\
&&+(-1)^{\alpha+\beta}T^{\alpha}_{~\delta}(1)\int_0^1 dt~U^{\delta}_{~\beta}(1,t)\big[\partial_{\theta}
\BS{\epsilon}^{\beta}_{~\gamma}
\big]\Big|_{\theta=0}(t)U^{\gamma}_{~\rho}(t,0)\nn\\
&&+(-1)^{\beta+\delta}\int_0^1 dt~U^{\alpha}_{~\beta}(1,t)\big[\partial_{\theta}
\BS{\epsilon}^{\beta}_{~\gamma}
\big]\Big|_{\theta=0}(t)U^{\gamma}_{~\delta}(t,0)T^{\delta}_{~\rho}(0).\nn\eea
The first two lines completely cancel, collecting the rest of the terms
\bea\delta_{\tiny{\textrm{chg of triv}}}U^{\alpha}_{~\gb}=\Big(-\epsilon^{\alpha}_{~\gd}(1)U^{\gd}_{~\gb}(1,0)+U^{\alpha}_{~\gd}(1,0)\epsilon^{\gd}_{~\gb}(0)\Big)
+(-1)^{\ga}\delta_BV^{\ga}_{~\gb}+(-1)^{\ga}T^{\ga}_{~\gd}(1)V^{\gd}_{~\gb}+(-1)^{\gd}V^{\ga}_{~\gd}T^{\gd}_{~\gb}(0)\nn\eea
where
\bea V^{\ga}_{~\gb}(1,0)=\int_0^1 dt~U^{\alpha}_{~\gd}(1,t)\big[(-1)^{\gd}\partial_{\theta}
\BS{\epsilon}^{\gd}_{~\gamma}
\big]\Big|_{\theta=0}(t)U^{\gamma}_{~\gb}(t,0).\nn\eea

\subsection{Examples from Courant Algebroid}
We use the setting of Ex.\ref{Courant}, we have a graded vector bundle
\bea &&T^*[2]T[1]M~~~~~~~~ (q_{\mu},p_{\mu})\nn\\
&&~~~~~~\downarrow\nn\\
&&~~~~T[1]M~~~~~~~(x^{\mu},v^{\mu})\nn\eea
This is a graded vector bundle in the sense that the transition function depends on the coordinate $v$, denote
\bea A^{\tilde\ga}_{\gb}=\frac{\partial\tilde x^{\tilde\ga}}{\partial x^{\gb}},\nn\eea
The transition function is
\bea \left[\tilde p_{\tilde\mu},~~\tilde q_{\tilde\mu}\right]
=\left[p_{\rho},~~q_{\rho}\right]\left[\begin{array}{cc}
   (A^{-1})^{\rho}_{\tilde\mu} & 0 \\
   -v^{\sigma}\partial_{\sigma}(A^{-1})^{\rho}_{\tilde\mu} & (A^{-1})^{\rho}_{\tilde\mu}
 \end{array}
\right].\label{transition}\eea
The homological vector field $Q=v^{\mu}\partial_{\mu}$ can be lifted
\bea \hat Q=Q+p_{\mu}\frac{\partial}{\partial q_{\mu}}.\nn\eea
It is not hard to check that such lifting is globally defined.
Namely $p_{\mu}\partial_{q_{\mu}}$ is a representation of $Q$ in the sense of Eq.\ref{representation}, the $T$ matrix
is
\bea T=\big[p_{\nu} ~~ q_{\nu}\big]
\left[\begin{array}{cc}
         0 & \delta_{\mu}^{\nu} \\
         0 & 0
       \end{array}\right]\left[\begin{array}{c}
                           \frac{\partial}{\partial p_{\mu}} \\
                           \frac{\partial}{\partial q_{\mu}}
                         \end{array}\right].\nn\eea
This constitutes the idea of adjoint representation up to homotopy of $T[1]M$, upon picking a connection. One can also check Eq.\ref{checked_later}
{\small\bea \Big[~\begin{array}{|cc|}
       A^{-1} & 0 \\
       -QA^{-1} & A^{-1}
     \end{array}\,,\begin{array}{|cc|}
                    1 & 0 \\
                    0 & -1
                  \end{array}\,
                  \begin{array}{|cc|}
                    0 & 1 \\
                    0 & 0
                  \end{array}~\Big]-\begin{array}{|cc|}
                    1 & 0 \\
                    0 & -1
                  \end{array}\,\begin{array}{|cc|}
                    QA^{-1} & 0 \\
                    0 & QA^{-1}\end{array}=0.\nn\eea}

Thus the Wilson-line is banally simple $U=\mathbf{1}$ in one trivialization patch. Under a change of trivialization, the rhs of Eq.\ref{triv_ind} had better be zero.
Writing $A=1+\epsilon$, then the $V$ term in that formula is
\bea V=\int_0^1 dt~\partial_{\theta}\left[\begin{array}{cc}
   -\BS{\ge} & 0 \\
   -\BS{Q\ge} & \BS{\ge}\end{array}\right]\bigg|_{\theta=0}
,\nn\eea
And we get for the three terms $(-1)^{\ga}\delta_BV^{\ga}_{~\gb}+(-1)^{\ga}T^{\ga}_{~\gd}(1)V^{\gd}_{~\gb}+(-1)^{\gd}V^{\ga}_{~\gd}T^{\gd}_{~\gb}(0)$
\bea
&&\int_0^1 dt~~\begin{array}{|cc|}
   1 & 0\\
   0 & -1\end{array}\,\begin{array}{|cc|}
   0 & 0\\
   -\partial_t(v^{\rho}\partial_{\rho}\ge) & 0\end{array}
   +\begin{array}{|cc|}
   1 & 0\\
   0 & -1\end{array}\,\begin{array}{|cc|}
   0 & 1\\
   0 & 0\end{array}\,\begin{array}{|cc|}
   0 & 0\\
   -\partial_t\ge & 0\end{array}
   +\begin{array}{|cc|}
   0 & 0\\
   -\partial_t\ge & 0\end{array}\,\begin{array}{|cc|}
   1 & 0\\
   0 & -1\end{array}\,\begin{array}{|cc|}
   0 & 1\\
   0 & 0\end{array}\nn\\
 &=&-\begin{array}{|cc|}
   \ge(1)-\ge(0) & 0\\
   Q\ge(0)-Q\ge(1) & \ge(1)-\ge(0)\end{array}\nn\eea
This term will cancel the first two terms of Eq.\ref{triv_ind}. We of course get the result '1 is invariant', but it is rather tricky to get the signs straight.

Next we check Eq.\ref{deform_not_simple} against the current example. It is obvious that the rhs of Eq.\ref{deform_not_simple} all vanishes for a Wilson-line in one trivialization patch. But the full Wilson-line is a product of transition functions
\bea U^{\ga_1}_{~\ga_n}=\Omega^{\ga_1}_{~\gb_1}(x_1)\Omega^{\gb_1}_{~\gb_2}(x_2)\cdots \Omega^{\gb_{n-1}}_{~\ga_n}(x_n),\nn\eea
but it suffices to consider just one transition function $U^{\ga}_{~\tilde\gb}=\Omega^{\ga}_{~\tilde\gb}(x)$. By a direct calculation, we have the deformation of the transition function
{\small\bea \partial_s\Omega=\begin{array}{|cc|}
                        (\bar v\cdot\partial)A^{-1} & 0 \\
                        -v^{\mu}\bar{v}^{\nu}\partial_{\mu}\partial_{\nu}A^{-1}-(d^y\bar v^{\mu})\partial_{\mu}A^{-1} & (\bar v\cdot\partial)A^{-1}\end{array}
=-d^y\begin{array}{|cc|}
                        0 & 0 \\
                        \bar v^{\mu}\partial_{\mu}A^{-1} & 0\end{array}+
\begin{array}{|cc|}
                        (\bar v\cdot\partial)A^{-1} & 0 \\
                        0 & (\bar v\cdot\partial)A^{-1}\end{array}\,,\label{temporary}\eea}
where we have omitted the index structure of $A^{-1}$, for which the reader may refer to Eq.\ref{transition}; and also in the second equalizer, we used the eom $v^{\mu}\partial_{\mu}=d^y$.

We should compare this to the general formula Eq.\ref{deform_trans}
\bea \partial_s\Omega^{\ga}_{~\tilde\gb}=d^y\iota_{\bar x}\Omega^{\ga}_{~\tilde\gb}+\iota_{\bar x}\Big(-T^{\ga}_{~\gd}\Omega^{\gd}_{~\tilde\gb}+(-1)^{\ga+\tilde\gc}\Omega^{\ga}_{~\tilde\gc}T^{\tilde\gc}_{~\tilde\gb}\Big).\nn\eea
The first term of this formula is obviously given by the first term of Eq.\ref{temporary}, while the second and third term can be worked out as
\bea&&\iota_{\bar x}\Big(-T^{\ga}_{~\gd}\Omega^{\gd}_{~\tilde\gb}+(-1)^{\ga+\tilde\gc}\Omega^{\ga}_{~\tilde\gc}T^{\tilde\gc}_{~\tilde\gb}\Big)\nn\\
&
=&\bar v^{\rho}\frac{\partial}{\partial v^{\rho}}\bigg(\,\begin{array}{|cc|}
                                                        0 & 1 \\
                                                        0 & 0
                                                      \end{array}\,
\begin{array}{|cc|}
                        A^{-1} & 0 \\
                        -(v\cdot\partial)A^{-1}& A^{-1}\end{array}+
\begin{array}{|cc|}
                                                        1 & 0 \\
                                                        0 & -1
                                                      \end{array}\,
\begin{array}{|cc|}
                        A^{-1} & 0 \\
                        -(v\cdot\partial)A^{-1}& A^{-1}\end{array}\,
\begin{array}{|cc|}
                                                        1 & 0 \\
                                                        0 & -1
                                                      \end{array}\,
\begin{array}{|cc|}
                                                        0 & 1 \\
                                                        0 & 0
                                                      \end{array}\,\bigg)\nn\\
&=&\begin{array}{|cc|}
                                                        (\bar v\cdot\partial)A^{-1} & 0 \\
                                                        0 & (\bar v\cdot\partial)A^{-1}
                                                      \end{array}\,,\nn\eea
which shows the total agreement with Eq.\ref{temporary}.

And finally, if $d^y=0$, then the first term of Eq.\ref{temporary} vanishes; and the second term also vanishes since we have shown immediately afterwards that it is in the image of the $T$ matrices, and see also the remark in Ex.\ref{information_loss}. Thus, we confirm that in the case $d^y=0$, the Wilson-line is invariant under the homotopies.


\end{document}